\documentclass{amsart}
\usepackage{amssymb, epsf, amscd}
\usepackage{xypic}                  
\allowdisplaybreaks
\swapnumbers
\newtheorem{theorem}[subsection]{Theorem}                
\newtheorem{mainthg0}[subsection]{Main Theorem $(g=0)$}
\newtheorem{mainthghi}[subsection]{Main Theorem $(g\ge0)$}
\newtheorem{proposition}[subsection]{Proposition}        
\newtheorem{lemma}[subsection]{Lemma}                    
\theoremstyle{definition}                      
\newtheorem{definition}[subsection]{Definition}          
\newtheorem{example}[subsection]{Example}                
\theoremstyle{remark}
\newtheorem{remark}[subsection]{Remark}                  
\numberwithin{equation}{section}
\newcommand{\ssref}[1]{Subsection~{\rm\ref{#1}}}
\newcommand{\firef}[1]{Figure~{\rm\ref{#1}}}
\newcommand{\thref}[1]{Theorem~\ref{#1}}
\newcommand{\prref}[1]{Proposition~\ref{#1}}

\newcommand{\deref}[1]{Definition~\ref{#1}}
\newcommand{\exref}[1]{Example~\ref{#1}}

\newcommand{\reref}[1]{Remark~\ref{#1}}

\newcommand{\fig}[1]
        {\raisebox{-0.5\height}%
                  {\epsfbox{legofigs/#1}}%
        }
\newcommand\captionfont[1]{\textsf{#1}}
\newcommand\lbb[1]{\label{#1} 
                   }                                    
\def\barF{$\bar{\text{F}}$}          
\def\barS{$\bar{\text{S}}$}          
\def\inv{${}^{-1}$}                  

\def\wti{\widetilde}
\def\ov{\overline}
\def\d{\partial}
\def\Max{{\mathrm{max}}}
\def\Min{{\mathrm{min}}}
\def\Sub{{\mathrm{Sub}}}
\def\edge{\rightsquigarrow}                    
\def\longedge{-\!\!\!-\!\!\!\rightsquigarrow}        
\def\isoto{\xrightarrow{\sim}}                 
\newcommand{\abs}[1]{\lvert#1\rvert}           

\def\i{{\mathrm{i}}}     
\def\Cset{\mathbb{C}}       
\def\Zset{\mathbb{Z}}       
\def\Qset{\mathbb{Q}}       
       
\def\al{\alpha}                         
\def\be{\beta}
\def\ga{\gamma}
\def\Ga{\Gamma}   
\def\de{\delta}

\def\ph{\varphi}

\def\Si{\Sigma}
\def\M{\mathcal{M}}                 
\def\Mt{\widetilde{\mathcal{M}}}
\def\Ct{\widetilde{\mathcal{C}}}
\def\C{\mathcal{C}}                  
\def\PGa{{\mathrm P}\Gamma}          
\DeclareMathOperator{\id}{id}     
\DeclareMathOperator{\HH}{H}            

\begin{document}

\title{On the Lego--Teichm\"uller Game}
\author{Bojko Bakalov}
\address{Department of Mathematics, MIT, Cambridge MA 02139}
\email{bakalov@math.mit.edu}
\thanks{The second author was supported in part by NSF Grant DMS9610201.}
\author{Alexander Kirillov, Jr.}
\address{Institute for Advanced Study, Princeton, NJ 08540, USA}
\email{kirillov@math.ias.edu}

\date{September 10, 1998}

\begin{abstract}
For a smooth oriented surface $\Si$, denote by $M(\Si)$ the set of all
ways to represent $\Si$ as a result of gluing together standard
spheres with holes (``the Lego game''). 
In this paper we give  a full set of simple moves and
relations which turn $M(\Si)$ into a connected and simply-connected
2-complex. Results of this kind were first obtained by Moore and
Seiberg, but their paper contains  serious gaps. Our proof is
based on a different approach and is much more  rigorous. 
\end{abstract}

\maketitle


\section{Introduction}\lbb{sintro}

Let $\Si$ be a smooth oriented surface, possibly with boundary. In
many cases---most importantly, for the study of the mapping class
group $\Ga(\Si)$ and for the construction of a modular functor---it is
convenient to represent $\Si$ as a result of gluing together several
simple pieces, which should be surfaces with a boundary. It is easy to
show that such a representation (we will call it a {\em
  parameterization}) is always possible if we allow these pieces to be
spheres with $\le3$ holes, or, more generally, spheres with $n$ holes.
For example, if we want to construct a modular functor, then it
suffices to define the vector spaces for each of these ``simple
pieces'', and then, since the behavior of the modular functor under
gluing is known, this defines uniquely the vector space which should
be assigned to $\Si$. From the point of view of the mapping class
groups, every parameterization defines a homeomorphism of the product
of the mapping class groups of the pieces into $\Ga(\Si)$. In
particular, in this way one can get a number of elements and relations
in $\Ga(\Si)$, and this can be used to get a full set of generators
and relations of $\Ga(\Si)$. This approach was first suggested by
Grothendieck (see below), who called it ``the Lego--Teichm\"uller
game''

In all of these applications, it is important to note that the same
surface $\Si$ can have many different parameterizations. Thus, it is
natural to ask the following questions. How can one describe
different ways of gluing ``standard pieces'' that give
parameterizations of the same surface $\Si$? Can we define some
``simple moves'' so that we can pass from a given parameterization to
any other by a sequence of these simple moves? And, finally, can one
describe all the relations between these simple moves, i.e.\ describe
when a sequence of simple moves applied to a parameterization yields
the same parameterization? 


These questions were studied in a series of pioneering papers of Moore
and Seiberg \cite{MS1, MS2}. These authors used spheres with 3 holes
(trinions) as their building blocks, and they gave a complete set of
simple moves and relations among them. However, their paper \cite{MS2}
has some serious flaws. First of all, they use the language of chiral
vertex operators, which is important for applications to conformal
field theory, but which is not really relevant for finding the set of
simple moves and relations, since this question is of purely 
topological nature. This lead them to miss some ``obvious'' axioms
which are automatically satisfied in any conforaml field theory.  What
is worse, their proof contains some gaps, the most serious of them
being a completely inadequate treatment of the case of surfaces of
higher genus with $n>1$ holes. The reason is that they used an
explicit presentation of the mapping class group $\Ga(\Si)$ by
generators and relations, found by Wajnryb \cite{Waj}, and such a
presentation for surfaces of higher genus was known only for surfaces
with $\le 1$ holes. (For surfaces with arbitrary number of holes, a
presentation of the mapping class groups by generators and relations
was recently found in \cite{Ge}; this presentation uses infinite
number of generators---all Dehn twists.)

In this paper, we give a reformulation and a rigorous proof of the
result of Moore and Seiberg, i.e., we construct a set of simple moves
and relations among them, which turn the set of all parameterizations
into a connected and simply-connected CW complex. To the best of our
knowledge, this has not been done accurately before; the only works in
this direction we are aware of are an unfinished and unpublished
manuscript by Kevin Walker, and the book \cite{T}, in which it is
proved that every modular tensor category  gives rise to a modular
functor.

We mostly use spheres with $n$ holes as our building blocks, which
allows us to simplify the arguments; however, we also give version of
the main theorem whihc only uses spheres with $\le 3$ holes. We do not
use explicit presentation of the mapping class groups $\Ga_{g,n}$ for
$g>0$ by generators and relations. Instead, we refer to the results of
Hatcher and Thurston \cite{HT} and their refinement by Harer \cite{H},
who solved a similar problem for the cut systems on $\Si$. Our
exposition is purely topological and requires no knowledge of modular
functors.

Our motivation for this work came form the conformal field theory and
modular functors. However, this work can also be useful for the study
of {\em Teichm\"uller tower\/}, introduced by A.~Grothendieck in his
famous {\em Esquisse d'un Programme\/} \cite{G}. This tower consists
of all stable algebraic curves of any genus $g\ge0$ with any number
$n\ge0$ of marked points (punctures) linked with the operation of
``gluing''.  The fundamental groupoid $T_{g,n}$ of the corresponding
moduli space is called the Teichm\"uller groupoid.  The most
fascinating thing is that the absolute Galois group
${\mathrm{Gal}}(\overline{\Qset}/{\Qset})$ acts on the profinite
completion $\{\widehat T_{g,n}\}$ of the tower of Teichm\"uller
groupoids, and this action is faithful---it is faithful already on
$\widehat T_{0,4}$ \cite{G} (see also \cite{S, LS} and references
therein).  Grothendieck states in \cite{G} as a very plausible
conjecture that the entire Teichm\"uller tower can be reconstructed
from the first two levels (i.e., the cases when $3g-3+n\le2$) via the
operation of ``gluing'', level $1$ gives a complete system of
generators, and level $2$ a complete system of relations:
\begin{quote}
\dots {\underline{la tour enti\`ere se reconstitue \`a partir des deux
premiers \'etages}}, en ce sens que via l'op\'eration fondamentale
de ``recollement'', l'\'etage $1$ fournit un syst\`eme complet de
g\'en\'erateurs, et l'\'etage $2$ un syst\`eme complet de relations.
\end{quote}

As Drinfeld says in \cite{D}, the above conjecture
``has been proved, apparently, in Appendix~B of the physics paper
\cite{MS2}''. As we already discussed, the approach of \cite{MS2} uses
heavily the explicit knowledge of the mapping class groups
and is not really rigorous. We believe that combining the results
of the present paper and of the unpublished manuscript \cite{BFM}, one can get
a proof of the Grothendieck conjecture. This will be discussed in
forthcoming papers.

\section{Extended surfaces and parameterizations}\lbb{sesurf}

\begin{definition}\lbb{dexsurf}
An {\em  extended surface} ({\em e-surface}, for short)
is a compact oriented smooth $2$-dimensional manifold
$\Si$, possibly with a boundary $\d\Si$,
with a marked point chosen on each boundary circle of $\Si$.
We will denote the set of boundary components
$\pi_0(\d\Si)$ by $A(\Si)$; in most cases, the elements of $A(\Si)$
will be labeled by Greek letters.

The {\em  genus\/} $g(\Si)$ of an e-surface $\Si$ is defined as
the genus of the surface without boundary, obtained from
$\Si$ by gluing a disk to every boundary circle.

A {\em morphism\/} of extended surfaces is an orientation-preserving
homeomorphism $\Si\isoto \Si'$ which maps marked points to marked points.
Every such morphism induces a bijection $A(\Si)\isoto A(\Si')$. 
%
\end{definition}


An example of an e-surface is shown in \firef{Fesf} below.

\begin{figure}[h]
\begin{equation*}
\fig{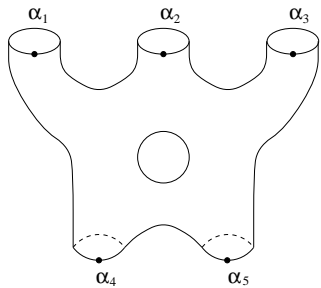}
\end{equation*}
\caption{\captionfont{An extended surface.}}\lbb{Fesf}
\end{figure}

\begin{definition}\lbb{dmcg} 
The {\em  mapping class group\/} 
$\Ga(\Si)$ of an e-surface $\Si$ is
the group of homotopy classes of morphisms 
$\Si\isoto\Si$. The 
{\em  pure mapping class group\/} $\PGa(\Si)$ is the subgroup
of $\Ga(\Si)$ of those morphisms that act trivially on 
the set of boundary components.
\end{definition}

\subsection{Standard sphere}\lbb{substsph}
For every $n\ge 0$, we define the {\em standard sphere\/} $S_{0,n}$ to be
the Riemann sphere $\ov{\Cset}$ with $n$ disks $|z-k|<1/3$ removed,
and with the marked points being $k-\i/3$ ($k=1,\dots,n$).  (Of
course, we could have replaced these $n$ disks by any other $n$
non-overlapping disks with centers on the real line and with marked
points in the lower half-plane---any two such spheres are
homeomorphic, and the homeomorphism can be chosen canonically up to
homotopy.) The standard sphere with 4 holes is shown in
\firef{Fstsphere}.  We will denote by $\Ga_{0,n}=\Ga(S_{0,n})$
(respectively, $\PGa_{0,n}=\PGa(S_{0,n})$) the mapping class group
(respectively, the pure mapping class group). 

\begin{figure}[h]
\begin{equation*}
\fig{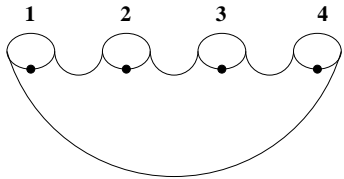}
\end{equation*}
\caption{\captionfont{A standard sphere (with $4$ holes).}}\lbb{Fstsphere}
\end{figure}

Note that the set
of boundary components of the standard sphere 
is naturally indexed by numbers $1, \dots, n$; we will use bold
numbers for denoting these boundary components:
$A(S_{0,n})=\{\mathbf 1, \dots, \mathbf n\}$.

Obviously, every connected e-surface of genus zero is homeomorphic to
exactly one of the standard spheres $S_{0,n}$, and the set of homotopy
classes of such homeomorphisms is a torsor over the mapping class
group $\Ga_{0,n}$.

\begin{definition}\lbb{dpar1}
Let $\Si$ be a connected e-surface of genus zero. A {\em parameterization
without cuts\/} of $\Si$ is a homotopy equivalence class of homeomorphisms
$\psi:\Si\isoto S_{0,n}$.
\end{definition}

The words ``without cuts'' are added in the definition because later
we will consider a more general notion of ``parameterization with
cuts''. 

We note that we have natural actions of the mapping class groups
$\Ga(\Si)$ and $\Ga_{0,n}$ on the set of all parameterizations of
$\Si$, given by $\ph (\psi)=\psi\circ \ph^{-1}, \ph\in \Ga(\Si)$ and
$\ph(\psi)=\ph\circ\psi, \ph\in \Ga_{0,n}$.  These actions are
transitive, so that the set of all parameterizations without cuts can
be identified (non canonically!) with $\Ga_{0,n}$. Note that every
parameterization defines an identification $A(\Si)\isoto\{\mathbf{1},
\dots, \mathbf{n}\}$ and thus, a natural order on the set of boundary
components $A(\Si)$.

\begin{definition}\lbb{dcut}
Let $\Si$ be an e-surface. A {\em cut system\/} $C$ on $\Si$ is a
finite collection of smooth simple closed non-intersecting curves on
$\Si$ (called {\em cuts}) such that each connected component of the
complement $\Si\setminus C$ is a surface of genus $0$.  (The cuts are
not to be oriented or ordered.) We will denote by $C(\Si)$ the set of
all cut systems on $\Si$ modulo isotopy.

A cut $c\in C$ is called {\em removable\/} if $C\setminus c$ is
again a cut system. A cut system is called {\em minimal\/}
if it contains no removable cuts.\footnote{Our notion of a minimal cut system
is exactly what is called a ``cut system'' in \cite{HT, H}.}
\end{definition}
Note that we could have defined a cut to be a simple closed curve with
one point marked on it. It is easy to see that this would have given
us the same set $C(\Si)$: given a cut $c$ and two points $p, p'$ on
$c$, there always exists an isotopy of $\Si$ which is different from
identity  only in a small neighborhood of $c$, maps $c$ onto itself
and $p$ to $p'$. 

Examples of cut systems and a minimal cut system on an e-surface
are shown in \firef{Fcut}.

\begin{figure}[h]
\begin{equation*}
\fig{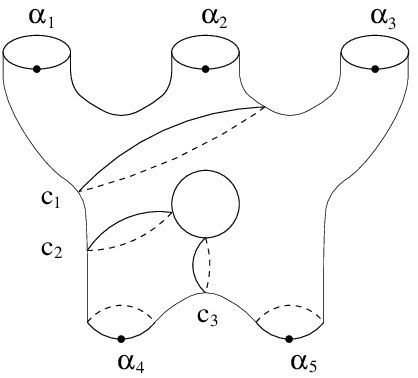}
\qquad\qquad
\fig{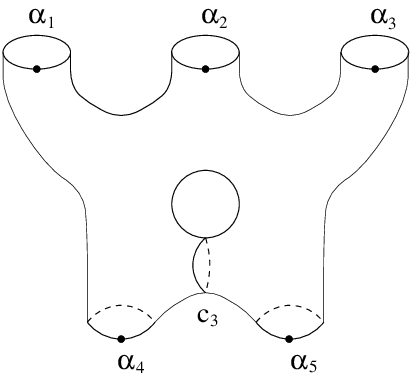}
\end{equation*}
\caption{\captionfont{Examples of cut systems and a minimal cut system 
(right).}}\lbb{Fcut}
\end{figure}

\begin{definition}\lbb{dpar2}
Let $\Si$ be an e-surface. A {\em parameterization} $P$ of $\Si$ is a
collection $(C, \{\psi_a\})$, where $C$ is a cut system on $\Si$, and
$\psi_a$ are parameterizations without cuts of the connected components
$\Si_a$ of $\Si\setminus C$, i.e.\ homotopy equivalence classes of
homeomorphisms $\psi_a:\Si_a\isoto S_{0, n_a}$ (see \deref{dpar1}).
\end{definition}

As before, we have a canonical action of $\Ga(\Si)$ on the set of all
parameterizations of $\Si$,  given by $\ph(C,
\{\psi_a\})=(\ph(C), \{\psi_a\circ \ph^{-1}\})$.

\section{Marked surfaces}\lbb{smsurf}
In this section, we will introduce some visual language for
representing the parameterizations of an e-surface $\Si$. 
Let us start with surfaces without cuts.

\begin{definition}\lbb{dmark1}
Let $S_{0,n}$ be the standard sphere with $n$ holes 
(see \ssref{substsph}), and let $m_0$ 
be the graph on it, shown in \firef{Fstmark} (for $n=4$).
This graph has a distinguished edge---the one which connects the
vertex $*$ with the boundary component $\mathbf 1$; in the figure,
this edge is marked by an arrow. 
We call $m_0$ the {\em standard marking without cuts\/} of
$S_{0,n}$. (For $n=0$, we let $m_0=\emptyset$.)

Let $\Si$ be an e-surface of genus zero. A {\em marking without
cuts\/} of $\Si$ is a graph $m$ on $\Si$ with one marked edge
such that  $m=\psi^{-1}(m_0)$ for some 
homeomorphism $\psi\colon\Si\isoto S_{0,n}$.  The graphs are
considered up to isotopy 
of $\Si$.
\end{definition}

\begin{figure}[h]
\begin{equation*}
\fig{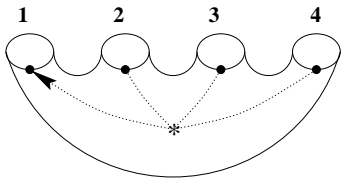}
\end{equation*}
\caption{\captionfont{Standard marking of the standard sphere 
(with $4$ holes).}}\lbb{Fstmark}
\end{figure}

Note that the ``free ends'', i.e., the vertices of the graph $m$ other
than $*$, coincide with the marked points on the boundary circles of
$\Si$.

\begin{proposition}\lbb{pmark1}
Let $\Si$ be an e-surface of genus zero. Then there is a bijection
between the set of all parameterizations without cuts of $\Si$
and the set of all markings without cuts of $\Si$,
given by $\psi\mapsto \psi^{-1}(m_0)$.
\end{proposition}
The proof of this proposition is elementary and is left to the
reader. We will denote either of the two sets of the proposition by
$M^\emptyset(\Si)$. 

Note that any marking without cuts of the surface $\Si$ defines a
bijection $A(\Si)\isoto \{\mathbf 1, \dots, \mathbf n\}$. In
particular, it defines an order on $A(\Si)$ and a {\em distinguished
boundary component}, corresponding to $\mathbf 1$.

Thus, these graphs provide a nice pictorial way of describing
parameterizations of e-surfaces. Similar to the constructions in the
previous section, we now define a more general notion of a marking with
cuts. 

\begin{definition}\lbb{dmark}
Let $\Si$ be an e-surface.  A {\em marking\/} $M$ of $\Si$ is a pair
$(C,m)$, where $C$ is a cut system on $\Si$ and $m$ is a graph on
$\Si$ with some distinguished edges such that it gives a marking
without cuts of each connected component of $\Si\setminus C$. We will
denote the set of all markings of a surface $\Si$ modulo isotopy by
$M(\Si)$. A {\em marked surface\/} ({\em m-surface\/}) is an e-surface
$\Si$ together with a marking $M$ on it.
\end{definition}

An example of an m-surface is shown in \firef{Fmsurf}.

\begin{figure}[h]
\begin{equation*}
\fig{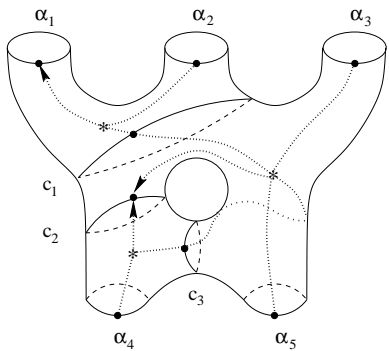}
\end{equation*}
\caption{\captionfont{A marked surface.}}\lbb{Fmsurf}
\end{figure}

 Note that, by definition, we have a canonical
``forgetting map'' $\pi\colon M(\Si)\to C(\Si)$ (recall that $C(\Si)$
is the set of all cut systems on $\Si$, see \deref{dcut}).

As before, the main reason for defining these markings is the following
result, which immediately follows from \prref{pmark1}. 

\begin{proposition}\lbb{pmark2}
Let $\Si$ be an e-surface. Then there is a bijection between the set of all 
parameterizations of $\Si$ and  the set of all markings of $\Si$. 
\end{proposition}

\subsection{Operations on markings}\lbb{subopmar}
 Rewriting the action of the mapping class group $\Ga(\Si)$ on the set
of all parameterizations of $\Si$, defined in \ssref{dpar1},
in terms of markings, we see that $\Ga(\Si)$ acts on $M(\Si)$ by
$\ph(C,m) = (\ph(C),\ph(m))$.  We also have the following obvious
operations:

\begin{description}
\item[Disjoint union]
$\amalg\colon M(\Si_1) \times M(\Si_2) \to M(\Si_1\amalg\Si_2)$.

\item[Gluing] If $\Si$ is an e-surface and $\al,\be\in A(\Si)$ is an
unordered pair with $\al\ne\be$, then we have a map
\begin{equation*}
\amalg_{\al,\be}\colon M(\Si) \to M(\Si'),
\end{equation*}
where $\Si' := \amalg_{\al,\be}(\Si)$ is obtained from 
$\Si$ by identifying the boundary
components $\al,\be$ so that the marked
points are glued to each other (this defines $\Si'$ uniquely up to
homotopy). The image of $\al$ and $\be$ is a cut on $\Si'$;
the marking on $\Si'$ is shown in \firef{Fglue}. If either of the
edges ending at $\al, \be$ was marked by an arrow, then we keep the
arrow after the gluing. 
\end{description}

\begin{figure}[h]
\begin{equation*}
\hspace{30pt}
\fig{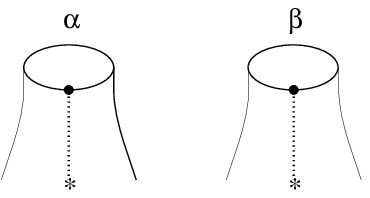}
\xrightarrow{ \;\; \sqcup_{\al,\be} \;\; }
\hspace{10pt}
\fig{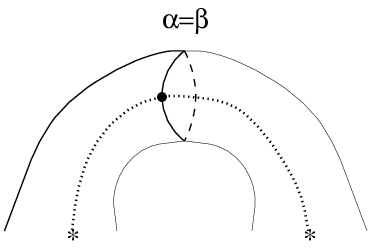}
\end{equation*}
\caption{\captionfont{Gluing of m-surfaces.}}\lbb{Fglue}
\end{figure}

The above two operations satisfy natural associativity properties, which
we do not list here (compare with \ssref{sprop}). Note also
that if $\Si_1$, $\Si_2$ are two e-surfaces, and $\al\in A(\Si_1)$,
$\be\in A(\Si_2)$, then we can define gluing
\begin{equation*}
\amalg_{\al,\be}\colon  M(\Si_1) \times M(\Si_2)
\to M(\Si_1\amalg_{\al,\be}\Si_2)
\end{equation*}
as the composition 
$M(\Si_1) \times M(\Si_2) \to M(\Si_1\amalg\Si_2) \to
M\bigl( \amalg_{\al,\be}(\Si_1\amalg\Si_2) \bigr) 
=:M(\Si_1\amalg_{\al,\be}\Si_2)$.

\subsection{Marking graphs}\lbb{smgraph}
Finally, note that for any m-surface $\Si$ and a marking $(C,m)\in M(\Si)$,
the graph $m$ has some
additional structure. Namely, the orientation of $\Si$ gives a natural
(counterclockwise) cyclic order on the set of germs of edges starting
at a given vertex.  Also, we have a distinguished set of $1$-valent vertices
(called ``free ends''), corresponding to the boundary components; these
vertices are in bijection with the set $A(\Si)$. It is also easy to
see that all ``internal edges''---i.e., edges that do not have a free
end---are in bijection with the cuts of $C$. 
We will always draw such graphs on the plane so that the cyclic
order on edges coincides with the counterclockwise order; as before,
we will mark the distinguished edges by arrows. 
For example, the surface in \firef{Fmsurf} gives the 
graph shown in  \firef{Fgraph}.

\begin{figure}[h]
\begin{equation*}
\fig{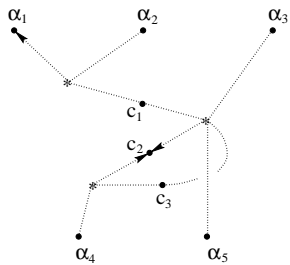}
\end{equation*}
\caption{\captionfont{A marking graph.}}\lbb{Fgraph}
\end{figure}

We also note that  Proposition~\ref{pmark2} easily implies the
following ``rigidity theorem''.

\begin{theorem}\lbb{trig} 
Let $\Si$ be a connected e-surface, $\al\in A(\Si)$, and 
$M=(C,m)$ --- a marking of  $\Si$. If
$\ph\in \Ga(\Si)$ is such that $\ph(M)=M, \ph(\al)=\al$, then $\ph=\id$.
\end{theorem}

\section{The complex $\M(\Si)$}\lbb{smsigma}

Let $\Si$ be an extended surface (see \deref{dexsurf}).  In
Subsections \ref{subg0mv}, \ref{sprop}, \ref{subg0rel},
\ref{subghimv}, \ref{subghirel} below, we will define a
$2$-dimensional CW complex $\M(\Si)$, which has the set $M(\Si)$ of
all markings of $\Si$ as the set of vertices.  The edges of $\M(\Si)$
will be directed; we call them {\em moves\/}.  It is convenient to
look at $\M(\Si)$ as a groupoid with objects---all vertices and
morphisms between two vertices---the set of homotopy classes of paths
on the edges of $\M(\Si)$ from the first vertex to the second one
(going along an edge in the direction opposite to its orientation is
allowed). We will use group notation writing a path composed of edges
$E_1, E_2, \dots$ as a product $E_1E_2\dotsm$, and we will write
$E^{-1}$ if the edge $E$ is traveled in the opposite direction. Then
the $2$-cells are interpreted as relations among the moves: we will
write $E_1\cdots E_k=\id$ if the closed loop formed by the  edges
$E_1,\dots, E_k$ is contractible in $\M(\Si)$; if we want to specify
the base point for the loop, we will write $E_1\dotsm E_k(M)=\id(M)$.
We will write $E:M\edge M'$ if the edge $E$ goes from $M$ to $M'$.

Our Main Theorems \ref{tmtg0} and \ref{tmtghi} 
state that the complex $\M(\Si)$ is connected and simply-connected.
Up to \ssref{subghimv}, the e-surface $\Si$ will be of genus $0$.
\subsection{Genus $0$ moves}\lbb{subg0mv}
We define the following {\em simple moves}:
\begin{description}

\item[Z-move]
Let $\Si$ be a connected e-surface of genus $0$ and $M=(\emptyset, m)$
a marking without cuts on $\Si$. Then we define the Z-move
$Z:(\emptyset, m)\edge (\emptyset, m')$, where $m'$ is the same graph
as $m$ but with a different distinguished edge, see \firef{Fzmv}.

\begin{figure}[h]
\begin{equation*}
\hspace{30pt}
\fig{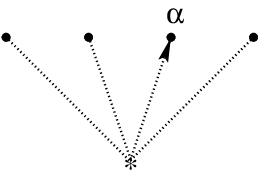}
\hspace{20pt}
\overset{Z}{\longedge}
\hspace{10pt}
\fig{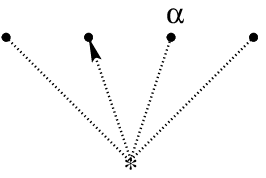}
\end{equation*}
\caption{\captionfont{Z-move (``rotation'').}}\lbb{Fzmv}
\end{figure}

\item[F-move]
Let $\Si$ be a connected  e-surface of genus $0$ and $(\{c\},m)\in M(\Si)$ 
be a marking with only one cut such that the edges ending at this cut
is  the distinguished (``first'') edge for
one of the components and the ``last'' edge for the other, as
illustrated in \firef{Ffmv} below.
Then we define  the F-move (``fusion'') $F_c:(\{c\},m)\edge(\emptyset,m')$, 
where the graph $m'$ is obtained from
 $m$ by contracting the edges ending at $c$, see \firef{Ffmv}.

\begin{figure}[h]
\begin{equation*}
\hspace{30pt}
\fig{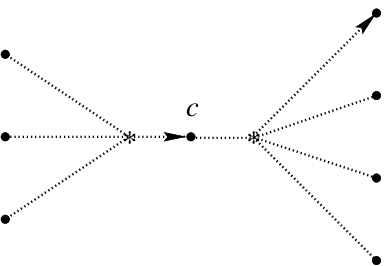}
\hspace{20pt}
\overset{F_c}{\longedge}
\hspace{10pt}
\fig{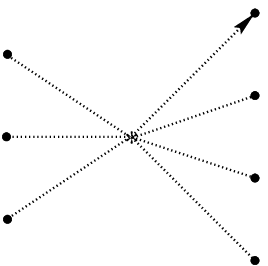}
\end{equation*}
\caption{\captionfont{F-move (``fusion'' or ``cut removal'').}}\lbb{Ffmv}
\end{figure}

\item[B-move]
Let $S_{0,3}$ be the standard sphere with three holes (a trinion), 
with no cuts and with the standard marking $m_0$,
shown in the left hand side of \firef{Fbmv}. 
We define the ``braiding'' move
$B_{\al,\be}$ by \firef{Fbmv}.

More generally, let $\Si$ be an e-surface and $\ph$ be a homeomorphism
$\ph\colon\Si\isoto S_{0,3}$. Then we define the move 
$B_{\al,\be} \colon \bigl( \emptyset, \ph^{-1}(m_0) \bigr) \edge 
\bigl( \emptyset, \ph^{-1}(B_{\al,\be}(m_0)) \bigr)$
in $\M(\Si)$.

\begin{figure}[h]
\begin{equation*}
\hspace{40pt}
\fig{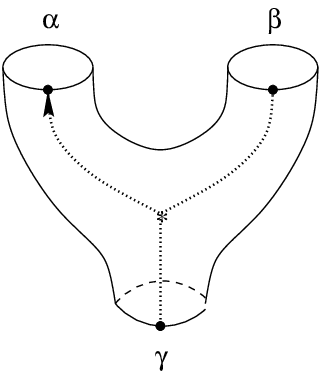}
\hspace{20pt}
\overset{B_{\al,\be}}{\longedge}
\hspace{10pt}
\fig{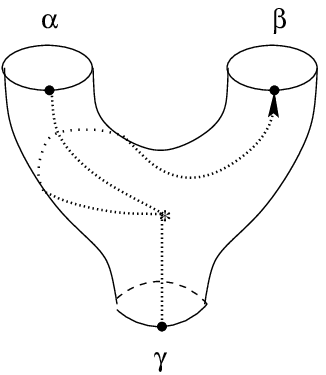}
\end{equation*}
\caption{\captionfont{B-move (``braiding'').}}\lbb{Fbmv}
\end{figure}

\end{description}

It easily follows from the definition that the set of moves is
invariant under the action of the mapping class groupoid: for every
edge $E:M_1\edge M_2, M_1, M_2\in M(\Si)$ and $\ph:\Si\to \Si'$, we
also have an edge $\ph(E):\ph(M_1)\edge\ph(M_2)$.

\subsection{Propagation of moves}\lbb{sprop}
We define the edges of the complex $\M(\Si)$ to be all that can be
obtained from the above defined simple moves by taking disjoint unions
and gluings modulo obvious equivalence relations. More precisely, we
define the set of edges of $\M(\Si)$ to be all edges that can be
obtained from the simple moves $Z, B, F$ by applying the two operations
listed below, modulo the equivalence relations below.

\textsc{Operations:} 
\begin{description}
\item[Disjoint union] If $\Si=\Si_1\sqcup\Si_2$, then for every edge
$E:M_1\edge M_1'$ in $\M(\Si_1)$, and every marking
$M_2\in M(\Si_2)$, we add an edge 
$E\sqcup\id_{M_2}:M_1\sqcup M_2\edge M'_1\sqcup M_2$ in $\M(\Si)$.

\item[Gluing] If $\Si_1=\sqcup_{\al,\be}\Si$,  then for every edge
$E:M\edge M'$ in $\M(\Si)$ we add an edge 
$\sqcup_{\al,\be}E: M_1 \edge M'_1$ in $\M(\Si_1)$, 
where $M_1=\sqcup_{\al,\be}M,
M'_1=\sqcup_{\al,\be}M'$, (cf.\ \firef{Fglue}).
\end{description}

\textsc{Equivalence relations:} 
\begin{description}
\item[Functoriality]
If $E,E'$ are edges in $\M(\Si_1)$, and $EE'$ is defined, then
\begin{equation*}
(E\sqcup\id_{M_2})(E'\sqcup\id_{M_2})
=(EE'\sqcup\id_{M_2}).
\end{equation*}
Similarly, we have
\begin{equation*}
\sqcup_{\al,\be}(E E') = (\sqcup_{\al,\be}E) (\sqcup_{\al,\be}E').
\end{equation*}

\item[Associativity 1]
For every edge $E$ in $\M(\Si_1)$ and markings $M_2\in M(\Si_2)$,
$M_3\in M(\Si_3)$, we have
\begin{equation*}
(E\sqcup\id_{M_2})\sqcup \id_{M_3}=E\sqcup\id_{M_2\sqcup M_3}.
\end{equation*}

\item[Associativity 2]
If $\al,\be,\ga,\de\in A(\Si)$ are four different boundary components
of $\Si$ and $E: M\edge M'$ is an edge, then 
\begin{equation*}
\sqcup_{\al,\be}(\sqcup_{\ga, \de} E)=\sqcup_{\ga, \de}(\sqcup_{\al,\be} E).
\end{equation*}

\item[Associativity 3] If $\Si=\Si_1\sqcup\Si_2, \al,\be\in A(\Si_1)$,
and $E:M\edge M'$ is an edge in $\M(\Si_1)$,  then 
\begin{equation*}
\sqcup_{\al,\be}(E\sqcup\id)=\sqcup_{\al,\be}(E)\sqcup\id.
\end{equation*}
\end{description}

\begin{remark}\lbb{redges}
Note that  an edge $E:M\edge M'$ just means
that we connect the points corresponding to $M$ and $M'$ in
$\M(\Si)$. It is {\bf not} a homeomorphism of surfaces. The relation
between the groupoid $\M(\Si)$ and the mapping class group $\Ga(\Si)$
can be described as follows: if $\ph\in \Ga(\Si)$ and  $M\in M(\Si)$
then we can ask if it is possible to connect $M$ with $\ph(M)$ by a
path in $\M(\Si)$. For example, the braiding move $B$  connects
$M$ with $b^{-1}(M)$ for a certain  $b\in \Ga(\Si)$ (``braiding''). 
One of the main results of the  next sections will be  that for every $\ph$,
there is a path in $\M(\Si)$, connecting $M$ with $\ph(M)$, and this
path is unique up to homotopy. However, as the example of the F-move
shows, we also have edges $E:M\edge M'$, where $M'$ cannot be obtained
from $M$ by the action of $\Ga(\Si)$.

To avoid confusion, we will denote elements of $\Ga(\Si)$ by lowercase
letters ($b,s,t,z,\dots$) and moves by uppercase letters ($B,
S,T,Z,\dots$). For the same reason, we use a different style of arrows
for edges.
\end{remark}

\begin{remark}\lbb{repaths}
When describing paths in the complex $\M(\Si)$, it is useful to note that
for every marking $M\in M(\Si)$ and any $\al,\be\in A(\Si)$,
there exists at most
one edge of the form $B_{\al,\be}$ originating from $M$: if we write
$B_{\al,\be}:M\edge M'$, then $M'$ is uniquely defined.  The same
applies to $Z$ and $F_c$. Thus, when describing a path in $\M(\Si)$ it
suffices to give the initial marking $M$ and a sequence of moves
$B^{\pm1},F^{\pm1},Z^{\pm1}$. 
This will define all the subsequent markings. However, to
assist the reader, in many cases we will make pictures of the
intermediate markings or at least of the corresponding graphs $m$.  
\end{remark}

\begin{example}[{\bf Generalized F-move}]\lbb{xgenf}
Let $\Si$ and $(\{c\},m)\in M(\Si)$ 
be as in the definition of the F-move, but with possibly
different distinguished edges for $m$. Let us fix an  order
of the connected components of $\Si\setminus c$, so that
$\Si=\Si_1\sqcup_c\Si_2$. We will call any composition of the form 
$ Z^a F_c (Z_1^k\sqcup\id)(\id\sqcup Z_2^l)$ a {\em generalized
  F-move}; for brevity, we will frequently denote it just by $F_c$. 
The Rotation axiom formulated below implies that up to homotopy, such
a composition is uniquely determined by the marking $M$ and by the
choice of the distinguished edge for the resulting marking
$F_c(M)$. Moreover, the Symmetry of F axiom along with the
commutativity of disjoint union, also formulated below, imply that if
we switch the roles of $\Si_1$ and $\Si_2$, then we get (up to homotopy)
the same generalized F-move. Thus, the homotopy class of the
generalized F-move is completely determined by the marking $M$ and by the
choice of the distinguished edge for the resulting marking
$F_c(M)$. 
\end{example}
\begin{example}[{\bf Generalized braiding}]\lbb{xgenbraid}
Let $\Si$ be a surface of genus zero, and let $m$ be a marking without
cuts of $\Si$. As discussed before, this defines an order on
the set of boundary components of $\Si$. Let us assume that we have a
presentation of $A(\Si)$ as a disjoint union, $A(\Si)=I_1\sqcup
I_2\sqcup I_3\sqcup I_4$, where the  order is given by $I_1<I_2<I_3<I_4$
(some of the $I_k$ may be empty).  Then we define the {\em generalized
braiding move\/} $B_{I_2, I_3}$ to be the product of simple moves
shown in \firef{Fgenbraid} below (note that we are using generalized
F-moves, see above).

\begin{figure}[h]
\begin{equation*}
\fig{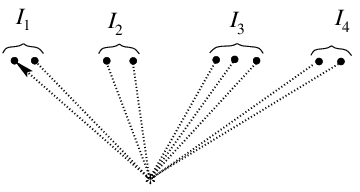}
\hspace{10pt}
\overset{F_{c_1}^{-1} F_{c_2}^{-1} F_{c_3}^{-1}}{\longedge}
\hspace{15pt}
\fig{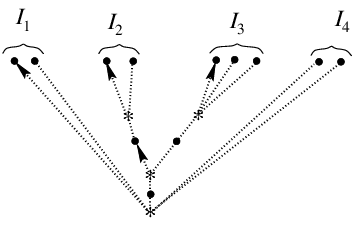}
\hspace{10pt}
\overset{B_{c_1, c_2}}{\longedge}
\end{equation*}
\begin{equation*}
\overset{B_{c_1, c_2}}{\longedge}
\hspace{15pt}
\fig{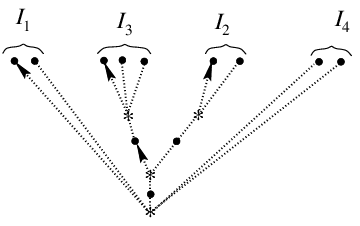}
\hspace{10pt}
\overset{F_{c_1} F_{c_2} F_{c_3}}{\longedge}
\hspace{15pt}
\fig{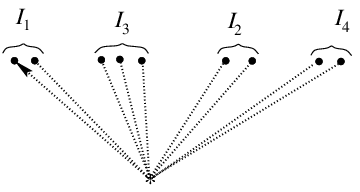}
\end{equation*}
\caption{\captionfont{Generalized braiding move.}}\lbb{Fgenbraid}
\end{figure}

It follows from the Cylinder axiom \eqref{cylax} below
that when $\Si$ is a three-holed sphere as in the definition of the
B-move, then  the B-move $B_{\al,\be}$
is homotopic in $\M(\Si)$ to the generalized braiding move 
$B_{\{\al\},\{\be\}}$.

We will use generalized moves  to simplify our formulas.
\end{example}

\subsection{Genus $0$ relations}\lbb{subg0rel}
Let us impose the following relations among the moves:
\begin{description}

\item[Rotation axiom]
If $\Si, M$ are as in the definition of the Z-move, then $Z^n=\id$,
where $n$ is the number of boundary components of $\Si$. 

\item[Commutativity of disjoint union]
If $\Si=\Si_1\sqcup\Si_2$ and $E_i$ is an edge in $\M(\Si_i)$ ($i=1,2$),
then in $\M(\Si)$
\begin{equation}\lbb{comax}
(E_1\sqcup\id)(\id\sqcup E_2)=(\id\sqcup E_2)(E_1\sqcup\id).
\end{equation}
We will denote either of these two products by $E_1\sqcup E_2$.

\item[Symmetry of  F-move]
Let $\Si, M$ be as in the definition of the F-move. 
Then $Z^{n_1-1} F_c=F_c(Z^{-1}\sqcup Z)$, where $n_1=|A(\Si_1)|$.

\item[Associativity of cuts] 
If $\Si$ is a connected surface  of genus zero, and $M=(C, m)\in M(\Si)$
is a marking with two cuts: $C=\{c_1, c_2\}$, then 
\begin{equation}\lbb{asax}
F_{c_1} F_{c_2} (M) = F_{c_2} F_{c_1} (M)
\end{equation}
(here $F$ denotes generalized F-moves).

\item[Cylinder axiom]
Let $S_{0,2}$ be a cylinder with boundary components $\al_0,\al_1$
and with the standard marking $M_0 = (\emptyset,m_0)$.
Let $\Si$ be an e-surface, $M=(C,m)\in M(\Si)$ be a marking, 
and $\al\in A(\Si)$ be a boundary component of $\Si$.
Then, for every move $E:M\edge M'$ in $\M(\Si)$
we require that the following square be contractible in 
$\M(\Si\sqcup_{\al,\al_1}S_{0,2})$:
\begin{equation}\lbb{cylax}
\begin{CD}
M \sqcup_{\al,\al_1} M_0 @>{E \sqcup_{\al,\al_1} \id}>> 
M' \sqcup_{\al,\al_1} M_0 \\
@V{F_\al}VV  @VV{F_\al}V \\
M @>>{\quad\;\; E \quad\;\;}> M'
\end{CD}\quad,
\end{equation}
where in the last line, we used the homeomorphism $\ph:
\Si\sqcup_{\al,\al_1}S_{0,2} \isoto \Si$, which is equal to identity
ouside of a neghborhood of $S_{0,2}$ and which maps $F_\al(M
\sqcup_{\al,\al_1} M_0)$ to $M$ (see \firef{Fcylax}). 

\begin{figure}[h]
\begin{equation*}
\fig{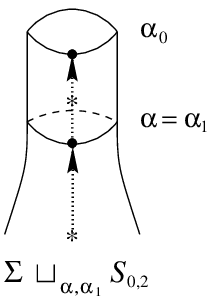}
\hspace{10pt} \xrightarrow{F_\al} \hspace{10pt}
\fig{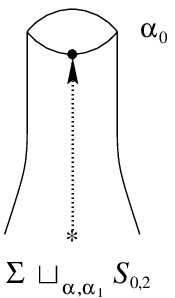}
\hspace{10pt} \xrightarrow{\;\;\ph\;\;} \hspace{10pt}
\fig{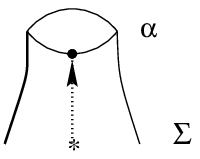}
\end{equation*}
\caption{\captionfont{Cylinder Axiom.}}\lbb{Fcylax}
\end{figure}

\item[Braiding axiom]
Let $\Si$ be an m-surface isomorphic to the sphere with $4$
holes, indexed by $\al$, $\be$, $\ga$, $\de$,  
and $M$  be a marking such that the graph $m$ is as shown in the left
hand side of \firef{Ftrax1}. 
Then
\begin{align}
\lbb{trax1}
B_{\al,\be\ga}(M) &= B_{\al,\ga} B_{\al,\be}(M),
\\
\lbb{trax2}
B_{\al\be,\ga}(M) &= B_{\al,\ga} B_{\be,\ga}(M).
\end{align}
For an illustration of Eq.~\eqref{trax1}, 
see Figure~\ref{Ftrax1}.
Note that all braidings involved are generalized braidings,
cf.\ \exref{xgenbraid}.

\begin{figure}[h]
\begin{equation*}
\hspace{20pt}
\fig{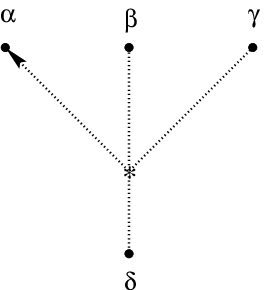}
\hspace{20pt}
\overset{B_{\al,\be}}{\longedge}
\hspace{10pt}
\fig{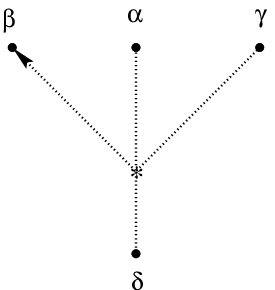}
\hspace{20pt}
\overset{B_{\al,\ga}}{\longedge}
\hspace{10pt}
\fig{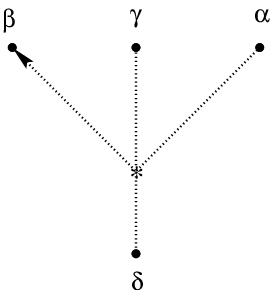}
\end{equation*}
\caption{\captionfont{Triangle axiom \eqref{trax1}.}}\lbb{Ftrax1}
\end{figure}

\item[Dehn twist axiom] Let $\Si$ be a sphere with $2$ holes, indexed
by $\al$, $\be$, and let $M=(\emptyset, m)$ be a marking without cuts
with the  distinguished vertex $\al$.
Then
\begin{equation}\lbb{axfort}
ZB_{\al,\be}(M) = B_{\be,\al}Z(M)
\end{equation}
(generalized braidings). This axiom is equivalent to the identity
$T_\al=T_\be$, where $T_\al$ is the Dehn twist
defined in \exref{xdehntw} below (see \firef{Fdehntw}).

\end{description}
\begin{proposition} \lbb{sendofdef}
All of the relations above make sense, i.e.\ they
describe closed paths in $\M(\Si)$.
\end{proposition} 
This proposition can be immediately verified explicitly.

We also add all relations that can be obtained from the above by taking
disjoint unions and results of gluing of relations: 
\begin{description}
\item[Propagation rules]
For every relation
$E=E'$ in $\M(\Si_1)$ we add  relations $E\sqcup\id=E'\sqcup\id$ in
$\M(\Si_1\sqcup\Si_2)$ and $\sqcup_{\al,\be} (E)=\sqcup_{\al,\be}(E')$
in $\M(\sqcup_{\al,\be}(\Si_1))$; compare with \ssref{sprop}.
\end{description}

This completes the definition of the complex $\M(\Si)$. Note that, by
definition, this complex is invariant under the action of the mapping
class groupoid: for every edge $E: M_1\edge M_2$ and $\ph:\Si\to \Si'$ we
also have an edge $\ph(E):\ph(M_1)\edge \ph(M_2)$. Similarly, for every
relation $E_1\dotsm E_n=\id$ we also have a relation
$\ph(E_1)\dotsm\ph(E_n)=\id$. 
 
Now we can formulate our main result for genus $0$ surfaces.
\begin{mainthg0}\lbb{tmtg0}
Let $\Si$ be an e-surface of genus $0$.
Then the above defined complex $\M(\Si)$ is connected
and simply-connected.
\end{mainthg0}
This theorem will be proved in Section~\ref{spfg0}. Here we give
several examples, which will play an important role later. 
\begin{example}[{\bf Associativity of cuts}]
If $\Si$ is a surface of genus zero, $M=(C,m)\in M(\Si)$ and $c_1,
c_2\in C$ are two of the cuts, then
\begin{equation}\lbb{asax1}
F_{c_1}F_{c_2}(M)=F_{c_2}F_{c_1}(M).
\end{equation}
Indeed, let us consider the connected components of
$\Si\setminus(C\setminus\{c_1,c_2\})$. If $c_1, c_2$ are in the same connected
component, then \eqref{asax1} follows from the Associativity axiom 
\eqref{asax} and
the Propagation rules. If $c_1$ and $c_2$ are in different connected
components, then \eqref{asax1} follows from the Commutativity of
disjoint union \eqref{comax}.
\end{example}
\begin{example}\lbb{xident0}
Let $\Si$ be an e-surface of genus zero and let $M$ be a marking with
one cut and with the marking graph shown in \firef{Fbf}. Then 
\begin{equation}
F_c B_{\{\al_1, \dots, c,\dots, \al_k\}, J}= 
B_{\{\al_1, \dots, I,\dots, \al_k\}, J} F_c.
\end{equation}
\begin{figure}[h]
\begin{equation*}
\fig{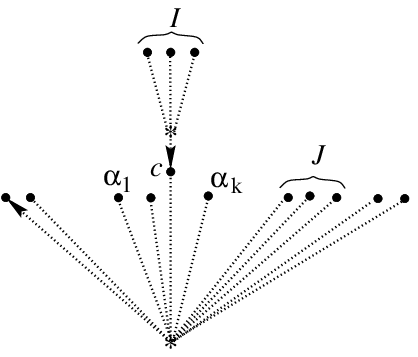}
\end{equation*}
\caption{}\lbb{Fbf}
\end{figure}
Indeed, this easily follows from the definition of the generalized
braiding and \eqref{asax1}. 
\end{example}
\begin{example}\lbb{xident1}
If we glue a disk to the hole $\be$ of the four-holed sphere
considered in the Braiding  axiom, we get 
$
B_{\al, \ga}=B_{\al, \ga}B_{\al,\emptyset},
$
which proves that $B_{\al,\emptyset}=\id$. We leave it to the reader
to write accurately the formal deduction of this from the above
axioms (this requires the use of the Associativity  axiom and the Commutativity
of disjoint union). Similarly, one can prove that
$B_{\emptyset, \al}=\id$, $B_{\emptyset,\emptyset}=\id$.  
\end{example}
\begin{example}[{\bf Generalized braiding}]\lbb{xident3}
Let $\al_1<\dots<\al_k<\be_1<\dots<\be_l$ be boundary components
of an m-surface $\Si$ of genus zero with no cuts. 
Then, applying repeatedly the Braiding axiom,
we get  
\begin{equation}\lbb{eqgenbr}
B_{ \{\al_1,\dots,\al_k\}, \{\be_1,\dots,\be_l\} }
= (B_{\al_1,\be_l}\dotsm B_{\al_k, \be_l})
\dotsm (B_{\al_1,\be_1}\dotsm B_{\al_k, \be_1}).
\end{equation}
(This argument also uses implicitly the Cylinder axiom and the
Propagation rules.)  

Hence, any path $B_{I_2, I_3}$ can be written as
a product of braiding moves of the form $B_{\al,\be}$ ($\al\in I_2,
\be\in I_3$). Note, however, that these $B_{\al,\be}$ are again
generalized braidings.
\end{example}
\begin{example}[{\bf Braid relations}]\lbb{xbrel}
In the setup of the Braiding axiom, one has the following relation:
\begin{equation}\lbb{eqbrel}
B_{\al,\be}B_{\al,\ga}B_{\be,\ga}=B_{\be,\ga}B_{\al,\ga}B_{\al,\be}.
\end{equation}
Indeed, by (\ref{trax1}, \ref{trax2}) this is equivalent to 
\begin{equation*}
B_{\al,\be}B_{\al\be, \ga}=B_{\be\al,\ga}B_{\al,\be},
\end{equation*}
which follows from the commutativity of disjoint union \eqref{comax}. 
By the Propagation rules,
it follows that \eqref{eqbrel} holds for a sphere with $n$ holes if
we choose as the basepoint a marking without cuts such that  
$\al<\be<\ga$.
\end{example}
\begin{example}[{\bf Dehn twist around a boundary component}]\lbb{xdehntw}
Let $\al$ be a boundary component of an m-surface $\Si$. 
Recall that the {\em Dehn
twist\/} $t_\al$ around $\al$ is the element of the mapping class
group $\Ga(\Si)$ which twists the boundary component $\al$ by 360
degrees counterclockwise. For any
marking $(C,m)$, we construct a path 
$T_\al: (C,m)\edge (C,t_\al^{-1}(m))$ on the edges of
$\M(\Si)$ as follows. 

First, we  define $T_\al$
for the standard cylinder $S_{0,2}$ with boundary components $\al,\be$
and the standard marking $(\emptyset, m_0)$ with a distinguished
vertex $\al$. 
Then let 
$T_\al :=B_{\al,\be}^{-1}Z$, where 
$B_{\al,\be} = B_{\{\al\},\{\be\}}$ 
is the generalized braiding from \exref{xgenbraid} 
(see \firef{Fdehntw}).

For an arbitrary marked surface $\Si$ and $\al\in A(\Si)$, define the
Dehn twist as the following composition: $Z^{-k}F_c T_\al F_c^{-1}
Z^k$, where $c$ is a small circle encompassing $\al$ and $k$ is such
that $Z^k(M)$ has $\al$ as the first boundary component. By the
Cylinder axiom, for $\Si$ being a cylinder this coincides with the
previously defined.

The  Dehn twist  axiom \eqref{axfort}
is equivalent to the identity $T_\al=T_\be$ for a cylinder with
boundary components $\al, \be$. 

\begin{figure}[h]
\begin{equation*}
\hspace{20pt}
\fig{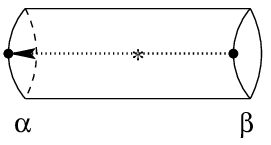}
\hspace{20pt}
\overset{T_\al}{\longedge}
\hspace{10pt}
\fig{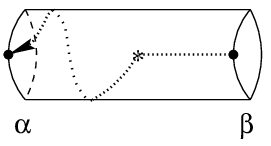}
\hspace{10pt}
=
\hspace{10pt}
\fig{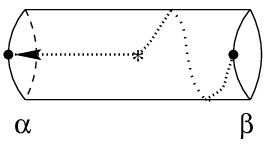}
\end{equation*}
\caption{\captionfont{Dehn twist ($T_\al=T_\be$).}}\lbb{Fdehntw}
\end{figure}

Similarly, if $\Si$ is a connected surface of genus $0$, $\al\in
A(\Si)$, and  $m$ is a
marking without cuts such that $A(\Si)=\al\sqcup I, \al<I$, then
\begin{equation}\lbb{eqdehn}
T_\al^{-1} = ZB_{\al,I}=B_{I,\al}Z^{-1}.
\end{equation}
This follows from the definitions and the  Dehn twist axiom.

By the commutativity of disjoint union, we have
\begin{equation}\lbb{dehncomm}
T_\al T_\be = T_\be T_\al
\end{equation}
for any $\al,\be\in A(\Si)$.
Also note that, by \exref{xident1}, we have $T_\al=\id$ for a sphere with one
hole $\al$, i.e.\ a disk. 
\end{example}
\begin{remark}\lbb{BTbt} One might ask why we chose $T_\al$ to connect
  $m$ with $t_{\al}^{-1} (m)$ rather then $t_\al(m)$. The reason is
  that if we recall that markings $m$ correspond to the homeomoprhisms
  $\psi:\Si\to S_{0,n}$ by $m=\psi^{-1} (m_0)$, then the marking
  $t_{\al}^{-1} (m)$ corresponds to the homeomorphism $t_i\circ \psi$,
  where $i$ is the index of the hole in $S_{0,n}$ corresponding to
  $\al$. Thus, the edge $T_\al$ connects $\psi$ with $t_i\circ \psi$.
  Similarly, the edge $B_{\al_i, \al_{i+1}}$ connects a
  homeomorphism $\psi$ with $b_i\circ\psi$, where $b_i\in
  \Ga(S_{0,n})$ is the braiding of $i$-th and $(i+1)$-st holes
  (cf. Proposition~\ref{pmcg0n}).
\end{remark}
\begin{example}[{\bf Dehn twist around a cut}]\lbb{xident2}
Let $\Si$ be a surface, and $M$ be a marking
containing a cut $c$. Define the moves $T'_c, T''_c$ by doing the
same construction as above on either of the sides of the
cut $c$. Then $T'_c=T''_c$. Indeed, it suffices to prove this when $\Si$ is
a cylinder, in which case it follows 
from the Cylinder axiom and the Dehn twist axiom that
$T'_c=T_\al$, $T''_c=T_\be$, $T_\al=T_\be$.
We will use the notation $T_c$ for both
$T_c',T''_c$. 
\end{example}
\begin{example}\lbb{xbalancing}
Let $\Si$ be a sphere with $3$ holes,
labeled by $\al,\be,\ga$, and let $M$ be a marking without cuts such
that it defines the  order $\al<\be<\ga$. Then we claim that
\begin{equation}\lbb{axbalancing}
T_\al(M) = T_\be T_\ga B_{\ga,\be}B_{\be,\ga}(M).
\end{equation}
Indeed, by \eqref{eqdehn} and the Braiding axiom, we have
\begin{equation*}
T_\al^{-1} = B_{\be\ga,\al} = B_{\be,\al}B_{\ga,\al}, \quad
T_\be^{-1} = B_{\be,\al}B_{\be,\ga}, \quad
T_\ga^{-1} = B_{\ga,\be}B_{\ga,\al},
\end{equation*}
which implies $T_\ga^{-1} T_\al T_\be^{-1} = B_{\ga,\be}B_{\be,\ga}$.
Now \eqref{axbalancing} follows from the commutativity
\eqref{dehncomm}.
\end{example}

\subsection{Higher genus moves}\lbb{subghimv}
Now let us consider e-surfaces $\Si$ of positive genus.
In this case, we need to add to the complex $\M(\Si)$
one more simple move and several more relations.
\begin{description}
\item[S-move]
 Let $S_{1,1}$ be a ``standard'' torus with one boundary
component and one cut, and with the  marking $M$  shown on the left
hand side of \firef{Fsmv}.  Then we add the edge
$S\colon M\edge M'$ where the marking $M'$ is shown on the right hand
 side of \firef{Fsmv}.

More generally, let $\Si$ be an e-surface and $\psi$ be a homeomorphism
$\psi\colon\Si\isoto S_{1,1}$. Then we add
the move $S\colon \psi^{-1}(M)\edge \psi^{-1}(M')$. 

\begin{figure}[h]
\begin{equation*}
\hspace{20pt}
\fig{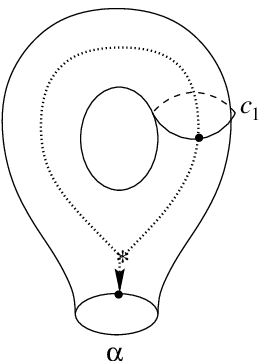}
\hspace{20pt}
\overset{S}{\longedge}
\hspace{10pt}
\fig{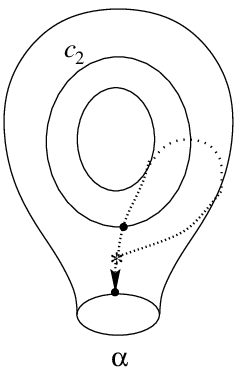}
\end{equation*}
\caption{\captionfont{S-move.}}\lbb{Fsmv}
\end{figure}

\end{description}

Of course, we also add all moves that can be obtained from
the above Z-, F-, B- and S-moves
by taking disjoint unions and gluing, as in \ssref{sprop}.

\begin{remark}\lbb{rsmv}
If $\Si$ is a surface of genus one with one hole,
we can identify the set of all markings with one cut on $\Si$
with the set of all homeomorphisms
$\psi\colon\Si\isoto S_{1,1}$ (see \thref{trig}).
Then the S-move connects  the marking
$\psi$ with $s\circ\psi$,
where $s\in\PGa(S_{1,1})$ acts on 
$\HH_1(\ov{S_{1,1}}) =\Zset c_1 \oplus\Zset c_2 \simeq\Zset^2$
by the matrix 
$\left(\begin{smallmatrix} 0 & 1 \\ -1 & 0 \end{smallmatrix}\right)$.
(Here $\ov{S_{1,1}}$ is the closed torus obtained by gluing a disk to
the boundary component of $S_{1,1}$.)
\end{remark}
\begin{example}[{\bf Generalized S-move}]\lbb{xgensmv}
Let $\Si$ be a torus with $n$ holes, and let $M=(\{c_1\}, m)$ be a
marking with one cut $c_1$, such that the graph $m$ is shown 
on the left hand side of
\firef{Fgensmv}. Then we define the {\em generalized
$S$-move\/} $S$ as the composition of moves shown in
\firef{Fgensmv}. 

\begin{figure}[h]
\begin{equation*}
\fig{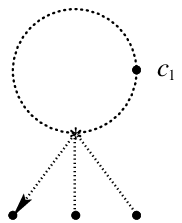}
\overset{F_c^{-1}}{\longedge}
\fig{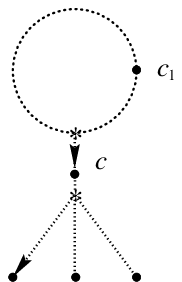}
\overset{S\sqcup_c\id}{\longedge}
\fig{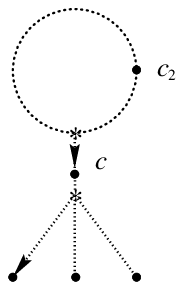}
\overset{F_c}{\longedge}
\fig{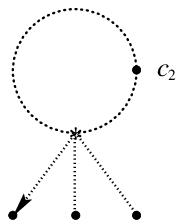}
\end{equation*}
\caption{\captionfont{Generalized S-move (for $n=3$).}}\lbb{Fgensmv}
\end{figure}

It can be shown that the cut $c$, and thus, the S-move, is uniquely
defined by $c_1$ and $m$ for $n\ge 1$.  As before, it follows from the
Cylinder axiom that for $n=1$, this generalized S-move coincides with
the one defined in \ssref{subghimv}. For $n=0$ there are two possible
choices for the cut $c$, and thus, there are two generalized S-moves
$S^{(1)}, S^{(2)}:M\edge M'$. It will follow from the relation
\eqref{g1n1r1} below that $S^{(1)}=S^{(2)}$. 
\end{example}

\subsection{Higher genus relations}\lbb{subghirel}
Besides those from \ssref{subg0rel},
we have the following additional relations:
\begin{description}
\item[Relations for $g=1,n=1$] Let $\Si$ be a marked torus with one hole
$\al$, isomorphic to the one shown on the left hand side of \firef{Fg1n1rel}.
For any marking $M=(\{c\}, m)$ with one cut, we let $T$ 
act on $M$ as the Dehn twist $T_c$ around $c$ (see \exref{xident2}).
Then we impose the following relations:
\begin{align}
S^2 &=Z^{-1} B_{\al,c_1},\lbb{g1n1r1}\\
(ST)^3&=S^2.\lbb{g1n1r2}
\end{align}
The left hand side of relation~\eqref{g1n1r1} is shown in 
\firef{Fg1n1rel}. For an illustration of \eqref{g1n1r2}, 
see Appendix~\ref{apg1n1r2}.

\begin{figure}[h]
\begin{equation*} 
\hspace{20pt}
\fig{smv1.eps}
\hspace{20pt}
\overset{S}{\longedge}
\hspace{5pt}
\fig{smv2.eps}
\hspace{20pt}
\overset{S}{\longedge}
\hspace{5pt}
\fig{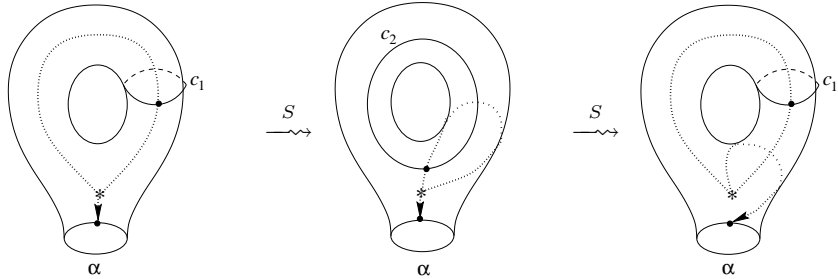}
\end{equation*}
\caption{\captionfont{The relation $S^2 = Z^{-1} B_{\al,c_1}$.}}
   \lbb{Fg1n1rel}
\end{figure}

\item[Relation for $g=1,n=2$]
Let $\Si$ be a marked torus with two holes $\al,\be$,
isomorphic to the one shown in \firef{Fg1n2}. 
As before, for any marking $M=(\{c\},m)$ of $\Si$ with one cut, we let $T$ act
on $M$ as the Dehn twist $T_c$ (see \exref{xident2}). Similarly, let 
$\wti T := T_{c +\be} = T_c T_\be B_{\be,c}B_{c,\be}$
(cf.\ Examples~\ref{xident2}, \ref{xbalancing}).
Then we have (note that we use generalized $S$-moves):
\begin{equation}
\lbb{g1n2}
 B_{\al,\be} F_{c_1} F_{c_2}^{-1} = S^{-1} {\wti T}^{-1} T S.
\end{equation}
The basepoint of the path \eqref{g1n2} is the marking shown in
\firef{Fg1n2} below. A detailed picture of the whole path is
presented in Appendix~\ref{apg1n2}.
\end{description}

\begin{figure}[h] 
\begin{equation*}
\fig{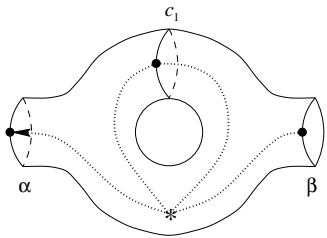}
\end{equation*}
\caption{\captionfont{A marked torus with two holes.}}\lbb{Fg1n2}
\end{figure}

Note that, by their construction, the above relations are
invariant under the action of the mapping class group.

It is not trivial that the relations (\ref{g1n1r2}, \ref{g1n2})
make sense, i.e.\ that they are indeed closed paths in $\M(\Si)$.
This is equivalent to checking that the corresponding identities
hold in the mapping class group $\Ga(\Si)$.
This is indeed so
(see, e.g., \cite{Bir, MS2}).
Another way of verifying that the relations
(\ref{g1n1r1}--\ref{g1n2}) make sense is to draw explicitly all the
marked surfaces and moves and to show that this is indeed a closed
loop.  We give such pictures in \firef{Fg1n1rel} and in
Appendices \ref{apg1n1r2} and \ref{apg1n2}.

\begin{example}\lbb{xtorus}
Let $\Si$ be a marked torus with one cut $c_1$
and one hole  $\al$ (see the left hand side of \firef{Fsmv}).
Then we have:
\begin{align}
\lbb{g1n1ri}
(ST)^3 &= S^2,
\\
\lbb{g1n1rii}
S^2 T &= T S^2,
\\
\lbb{g1n1riii}
S^4 &= T_\al^{-1}.
\end{align}
Indeed, \eqref{g1n1ri} is exactly \eqref{g1n1r2}.
Equation~\eqref{g1n1rii} follows from \eqref{g1n1r1}, the Cylinder
axiom, and the commutativity of disjoint union, and \eqref{g1n1riii}
easily follows from \eqref{g1n1r1} and \eqref{eqdehn}.

In particular, this implies that the elements $t, s\in \Ga_{1,1}$
(cf.\ \reref{redges})
satisfy the relations (\ref{g1n1ri}--\ref{g1n1riii}). In fact, it
is known that these are the defining relations of the group
$\Ga_{1,1}$ (see \cite{Bir}). 
\end{example}

Now we can formulate our main result for arbitrary genus.
\begin{mainthghi}\lbb{tmtghi}
Let $\Si$ be an e-surface.  Let $\M(\Si)$ be the above defined complex
with vertices: all markings of $\Si$, edges obtained from the Z-, F-, B-,
and S-moves by disjoint unions and gluing, and $2$-cells given by the
relations in Subsections {\rm\ref{subg0rel}}, {\rm\ref{subghirel}}.  Then
$\M(\Si)$ is connected and simply-connected.
\end{mainthghi}
This theorem will be proved in Section~\ref{spfghi}.

\section{The complex $\M^\Max(\Si)$}\lbb{smmaxsi}

In this section, we formulate and prove a version of the Main
Theorem~\ref{tmtghi} in which one uses only spheres with $\le3$
holes.
\footnote{This is referred to as ``a small lego box'' by
  Grothendieck \cite{G}.}

Clearly, any surface $\Si$ can be cut into a union of 
spheres with $\le3$ holes, i.e.\ trinions, cylinders, disks, and spheres.
Let us call such a cut system on $\Si$ {\em maximal}.
We will define a $2$-dimensional CW complex $\M^\Max(\Si)$
with vertices: all markings $M=(C,m)\in M(\Si)$ such that the 
cut system $C$ is maximal (we call such $M$ maximal).
The edges and the $2$-cells of $\M^\Max(\Si)$ 
are described below.

\textsc{Edges:}
\begin{description}
\item[Z-move] defined in the same way as in
  Subsection~\ref{subg0mv}, but for spheres with $\le 3$ holes. 

\item[F-move] $F_c:M\edge M'$ -- defined as in
   Subsection~\ref{subg0mv}, but for $M, M'\in\M^\Max(\Si)$. It is
   easy to see that this happens iff $\Si=\Si_1\sqcup_c\Si_2$ with
   only one cut $c$ and one of $\Si_1, \Si_2$ is either a cylinder or a disk.
If $\Si=\Si_1\sqcup_c\Si_2$ and both $\Si_1$ and $\Si_2$ are trinions,
then $F_c$ is not defined in $\M^\Max(\Si)$, because the result
would be a sphere with $4$ holes without cuts.

\item[A-move] if $\Si$ is a sphere with
$4$ holes, and $M$ is a marking with one cut shown in the left-hand
side of \firef{Fasmv} below, then we define the A-move $A:M\edge M'$
as in \firef{Fasmv}.

\item[B-move] defined in the same way as in Subsection~\ref{subg0mv}.

\item[S-move] defined in the same way as in Subsection~\ref{subghimv}.
\end{description}

As before, we also add all edges which can be obtained from these ones
by disjoint unions and gluing modulo the equivalence relations of
Subsection~\ref{sprop}. 

Thus, the edges of the complex $\M^\Max(\Si)$ are those edges of
$\M(\Si)$ which have both endpoints in $M^\Max(\Si)$, plus the new
A-moves. Note that for every A-move $A_{c',c}:M\edge M'$, 
the same markings can
be connected in $\M(\Si)$ by the path $F^{-1}_{c'} F_c(\Si)$.

\begin{figure}[h]
\begin{equation*}
\hspace{20pt}
\fig{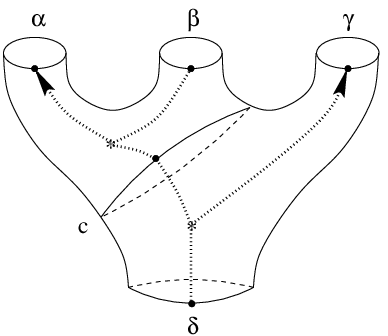}
\hspace{20pt}
\overset{A_{c',c}}{\longedge}
\hspace{10pt}
\fig{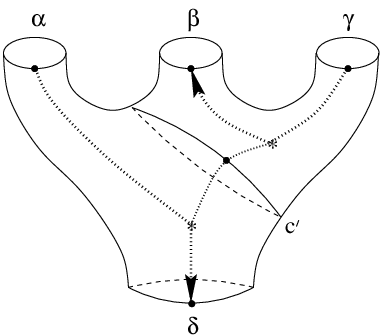}
\end{equation*}
  \caption{\captionfont{A-move (``associativity constraint'').}}\lbb{Fasmv}
\end{figure}

Note that the generalized braiding moves can not be defined in
$\M^\Max(\Si)$, with the exception of the Dehn twist $T_c$: its
definition can be repeated in $\M^\Max$.

\textsc{ Relations:}
\begin{description}
\item[Weak associativity of cuts]
if $\Si$ is a surface of genus zero, and $M\in \M^\Max(\Si)$ is a
marking with two cuts $c_1,c_2$ such that
$\Si=\Si_1\sqcup_{c_1}\Si_2\sqcup_{c_2}\Si_3$, and $\Si_2$ is a
cylinder, then $F_{c_1}=F_{c_2}$. 

\item[Symmetry of  F] same as in
Subsection~\ref{subg0rel}, but one of  $\Si_1, \Si_2$ has at most 2
holes, and the other at most 3. 

\item[Rotation axiom] same as in
Subsection~\ref{subg0rel}, but for $n\le 3$.

\item[Commutativity of disjoint union] same as in
Subsection~\ref{subg0rel}.

\item[Cylinder axiom] same as in
Subsection~\ref{subg0rel}. 

\item[Pentagon axiom](see Appendix~\ref{appenthex}).

\item[Two Hexagon axioms](see Appendix~\ref{appenthex}).

\item[Self-duality of associativity] if $\Si, M$ are as in definition of
the A-move, then $A^2(M)=\id(M)$. 

\item[Triangle axiom](see Appendix~\ref{appenthex}).

\item[Dehn twist axiom] same as in
Subsection~\ref{subg0rel}.

\item[Relations for $g=1, n=1$] same as in
Subsection~\ref{subghirel}.

\item [Relation for $g=1, n=2$](see Appendix~\ref{apg1n2}).
\end{description}

As before, we also add all the relations that can be obtained from
these ones by disjoint unions and gluing. This completes the
definition of the complex $\M^\Max(\Si)$.

Note that the Hexagon axioms are essentially the Braiding axioms
\eqref{trax1}, \eqref{trax2}, only rewritten so that they start at a
maximal marking and instead of products of the form $F^{-1}_{c'} F_c$
we used the  A-moves. 
Indeed, by the definition of generalized braiding (\exref{xgenbraid}),
the left hand side of \eqref{trax1} is $F_{c_2} B_{\al,c_2} F_{c_2}^{-1}$,
while the right hand side is 
$F_{c_3} B_{\al,\ga} F_{c_3}^{-1} F_{c_1} B_{\al,\be} F_{c_1}^{-1}$.
Therefore, \eqref{trax1} can be rewritten as 
\begin{equation}\lbb{eqhex}
B_{\al,c_2} A_{c_2,c_1} = A_{c_2,c_3} B_{\al,\ga} A_{c_3,c_1} B_{\al,\be},
\end{equation}
which is  the Hexagon relation.

The same can be said about the relation for $g=1,n=2$, 
see Appendix~\ref{apg1n2}.
Thus, the only new relations are the Pentagon and Triangle axioms. 

\begin{theorem}
For any e-surface $\Si$, the complex $\M^\Max(\Si)$ defined above is
connected and simply-connected.
\end{theorem}

\begin{proof}
It suffices to show that the complexes $\M^\Max(\Si)$ and $\M(\Si)$
are homotopically equivalent, after which the result follows from
Theorem~\ref{tmtghi}. To show the equivalence, we introduce the
following notion. Let $M,M'\in M(\Si)$. We say that $M'$ is a
subdivision of $M$ if $M$ can be obtained by applying to $M'$ a
sequence of F-moves (not $F^{-1}$!) and Z-moves. We will write
$M\subset M'$. 

Now, for a given $M\in M(\Si)$, denote 
\begin{equation*}
\Sub(M)=\{M'\in M^\Max(\Si)|M\subset M'\}.
\end{equation*}

\begin{lemma} 
\begin{enumerate}
\item Every two markings $M', M''\in \Sub(M)$ can be connected by a
  path in $\Sub(M)$ consisting of a sequence of the F-, Z-, and
  A-moves and their inverses.

\item Every loop in $\Sub(M)$, composed of F-, Z-, and A-moves  and
their inverses, is contractible in $\M^\Max(\Si)$. 
\end{enumerate}
\end{lemma}

The proof of this lemma is left to the reader. Obviously, it suffices
to consider the case $\Si=S_{0,n}, M=(\emptyset, m_0)$, in which case
it is essentially a version of  MacLane's coherence theorem.

Now, let us choose for every $M\in M(\Si)$ one element $\tau(M)\in
\Sub(M)$ (``maximal subdivision'') in such a way that we do not add
any new cuts to components which already are spheres with $\le 3$
holes. Then one easily sees that the map $\tau: M(\Si)\to
M^\Max(\Si)$ can be extended to a map of CW complexes. Indeed, it is
obvious how this map is defined on the edges of B and S type. As for
the F-edge, let us define for $F_c:M_1\sqcup_c M_2\edge M$ its image
$\tau(F_c)$ as any path $\tau(M_1) \sqcup_c \tau(M_2)\edge \tau(M)$,
composed of A-, Z-, and F-moves in $\M^\Max$; by the lemma above, such
a path is unique up to homotopy. Similarly we define $\tau(Z)$. It is
immediate to see that $\tau$ respects all the relations in $\M$.

Conversely, we have an obvious embedding $\M^\Max\subset \M$. It is
immediately verified that the composition of this embedding with $\tau$ is
an auto-equivalence of $\M^\Max$. Thus, every loop in $\M^\Max$ is
homotopic to a loop of the form $\tau(l)$ for some closed loop in
$\M$. But every loop $l$  in $\M$ is contractible; thus, the same holds for 
$\tau(l)$. 
\end{proof}

\section{Proof of the main theorem for genus $0$}\lbb{spfg0}

In this section we will prove the Main Theorem~\ref{tmtg0}
for genus $0$: that for any extended surface $\Si$ of genus $0$
the complex $\M(\Si)$, defined in 
Subsections~\ref{subg0mv}, \ref{sprop}, \ref{subg0rel},
is connected and simply-connected.

\subsection{Outline of the proof}\lbb{soutpf}
Let us start by slightly modifying the complex $\M(\Si)$. Namely, let
us add to it edges corresponding to each of the paths $B_{I_1, I_2}$
defined in \exref{xgenbraid} and the Dehn twists $T_\al$ (see examples
\ref{xdehntw}, \ref{xident2}). 

For each of these new moves, we add the expression for it as a product
of simple moves $Z,F,B$ as a new relation.  As before, we also add all
edges and relations that can be obtained from these ones by disjoint
union and gluing. Let us call the complex obtained in this way
$\Mt(\Si)$. Obviously, if $\Mt(\Si)$ is connected and
simply-connected, then so is $\M(\Si)$.

The proof is based on extending the forgetting map $\pi: M(\Si)\to
C(\Si)$ (see \ssref{dmark}) to a map of CW complexes $\pi:
\Mt(\Si)\to \C(\Si)$, and showing that both the base and any fiber are
connected and simply-connected. More precisely, we will use the
following proposition, whose easy proof is left to the reader. For
future use, we formulate it in a slightly more general form than we
need now.

\begin{proposition}\lbb{pmpic}
Let $\M$, $\C$ be $2$-dimensional CW complexes {\rm(}with directed edges{\rm)},
and let $\pi\colon \M^{[1]} \to \C^{[1]}$
be a map of their $1$-skeletons, which is surjective both
on vertices and on edges. Suppose that the following conditions are
satisfied{\rm:}
\begin{enumerate}
\item
$\C$ is connected and simply-connected.

\item
For every vertex $C\in \C$, $\pi^{-1}(C)$ is connected and
simply-connected  in $\M$ {\rm(}that is, every closed loop $l$ which lies
completely in $\pi^{-1}(C)$ is contractible in $\M${\rm)}.

\item
Let $C_1 \overset{e}{\edge} C_2$ be an edge in $\C$, and let $M'_1
\overset{e'}{\edge} M'_2$ and $M''_1 \overset{e''}{\edge} M''_2$ be
two its liftings to $\M$. Then one can choose paths $M'_1
\overset{e_1}{\edge} M''_1$ in $\pi^{-1}(C_1)$ and $M'_2
\overset{e_2}{\edge} M''_2$ in $\pi^{-1}(C_2)$ such that the square
\begin{equation*}
\begin{CD}
M'_1 @>{e'}>> M'_2 \\
@V{e_1}VV  @VV{e_2}V \\
M''_1 @>>{e''}> M''_2
\end{CD}
\end{equation*}
is contractible in $\M$.

\item
For every $2$-cell $X$ in $\C$, its boundary $\d X$ can be lifted
to a contractible loop in $\M$.
\end{enumerate}
Then the complex $\M$ is connected and simply-connected.
\end{proposition}

\subsection{The complex $\C(\Si)$} \lbb{subcg0} 
%
The set of vertices of $\C(\Si)$ is the set $C(\Si)$ of all cut systems
on $\Si$.
The (directed) edges of $\C(\Si)$ will correspond to
the following
\begin{description}
\item[\barF-move] 
Let $\Si$ be an e-surface of genus zero, and let 
$C\in C(\Si)$ be a cut system on $\Si$, consisting of a single cut:
$C=\{c\}$.   Then we define  the
\barF-move $\bar F_c \colon C \edge \emptyset$, which removes $c$.

\end{description}

As before, we also add all the moves that can obtained from the 
\barF-move above by disjoint unions and gluing subject to the 
obvious associativity relations as in  \ssref{sprop}. In particular,
for any $C\in C(\Si)$ ($\Si$ not necessarily of genus zero) and
a removable cut $c\in C$, we have a move $\bar F_c:C\edge C\setminus
\{c\}$. 

Let us impose the following relations for these moves:

\begin{description}
\item[Associativity of cut removal]
Let $c_1,c_2 \in C$, $c_1\ne c_2$. Then 
\begin{equation*}
\bar F_{c_1}\bar F_{c_2}(C)=\bar F_{c_1}\bar F_{c_2}(C).
\end{equation*}
\end{description}

We add these relations, as well as all relations obtained by taking disjoint
unions and gluing (cf.\ \ssref{sendofdef}),
as $2$-cells of the complex $\C(\Si)$.

By construction, there is a canonical map of CW complexes
$\pi\colon\Mt(\Si)\to\C(\Si)$, which extends the forgetting map
$\pi\colon M(\Si)\to C(\Si): (C,m)\mapsto C$. Namely, we define $\pi$
on edges by $\pi(F) = \bar F$, $\pi(B) = \id$, $\pi(Z)=\id$. 

\begin{theorem}[$g=0$]\lbb{tcsig0}
The above complex $\C(\Si)$ is connected and simply-connected.
\end{theorem}
\begin{proof}
It is easy to see from the Associativity axiom that every product
$\bar F_{c_1}\bar F_{c_2}^{-1}$ can be replaced by either $\bar
F_{c_2}^{-1}\bar F_{c_1}$ or by identity.  Thus, every loop can be
deformed to one of the form $\bar F^{-1}\dotsm \bar F^{-1}\bar F\dotsm
\bar F$.  On the other hand, every cut system can be connected to the
empty one using \barF-moves (this is where we need that $\Si$ is of
genus zero!). Thus, it suffices to consider only loops starting at the
empty cut system. But every loop of the form $\bar F^{-1}\dotsm \bar
F^{-1}\bar F\dotsm \bar F$ starting at the empty cut system must be
homotopic to identity.
\end{proof}

\subsection{Simply-connectedness of the fiber}\lbb{subscfiber}
Let $C$ be a vertex of $\C(\Si)$, i.e.\ a cut system on $\Si$.  Denote
by $\{\Si_a\}$ the set of connected components of $\Si\setminus C$.
Then $\pi^{-1}(C) \subset M(\Si)$ can be canonically identified with
the product $\prod_a M^\emptyset(\Si_a)$, where $M^\emptyset(\Si_a)$
is the set of all markings without cuts of $\Si_a$, (cf.\ 
\deref{dmark1}).  Thus, to check assumption 2 of
Proposition~\ref{pmpic}, it is enough to check that every
$\Mt^\emptyset(\Si )$, where $\Si $ is a sphere with $n $ holes, is
connected and simply-connected.  (Here $\Mt^\emptyset(\Si)$ is the
subcomplex of $\Mt(\Si )$ with vertices $M^{\emptyset}(\Si )$, and
edges given by Z-moves and the generalized B-moves.)

By \prref{pmark1}, the set
$M^\emptyset(\Si )$ is in bijection with the mapping class group
$\Ga_{0,n }= \Ga(S_{0,n })$. Let us consider the following elements of
$\Ga_{0,n }$: 
\begin{align*}
t_i,& i=1,\dots, n&&:\text{Dehn twist around $i$-th puncture}\\
b_i,& i=1,\dots, n-1&&:\text{Braiding of $i$-th, $(i+1)$-st punctures}\\
z&&&:\text{Rotation, i.e.\ a homeomorphism which acts  on the set}\\
 &&&\quad \text{of boundary components by $\mathbf i\mapsto\mathbf{i+1},
   \mathbf n\mapsto \mathbf 1$}\\
 &&&\quad\text{and preserves the real axis.}
\end{align*}
\begin{proposition}\lbb{pmcg0n}
The group $\Ga_{0,n }$ is 
generated by elements 
$b_i$, $i=1,\dots, n-1$, 
$t_i$, $i=1,\dots, n$, and $z$
with the following defining relations

\begin{align}
\lbb{ga0n1}
        &b_ib_j=b_jb_i, \qquad\qquad\qquad \abs{i-j}>1,\\
\lbb{ga0n2}
        &b_ib_{i+1}b_i=b_{i+1}b_ib_{i+1},\\
\lbb{ga0n3}
        &b_i t_j=t_jb_i, \qquad\qquad\qquad \abs{i-j}>1, \; i=j+1,\\
\lbb{ga0n4}
        &b_i^{\pm 1}t_i=t_{i+1}b_i^{\pm 1},\\
\lbb{ga0n5}
        &t_i t_j=t_j t_i,\\
\lbb{ga0n6}
        &z^n=1,\\
\lbb{ga0n7}
        &b_1\dots b_{n-1}t_n=z,\\
\lbb{ga0n8}
        &zt_n=t_1z.
\end{align}

\end{proposition}

This proposition is known (see, e.g., \cite{MS2}, where it is
formulated in a somewhat different form), so we skip the proof.

\begin{remark}\lbb{remcg0n}
Denote the boundary components of $\Si$ by $\al_1,\dots,\al_n$.
Let $\psi$ be a homeomorphism $\Si\isoto S_{0,n }$ which
induces the  order $\al_1<\dots<\al_n$. As was noted
before, such a homeomorphims can be viewed as an element of  
$M^\emptyset(\Si)$. Then in
$\Mt^\emptyset(\Si)$ we have the edges

\begin{equation}\lbb{btmoves}
\begin{aligned}
B_{\al_i,\al_{i+1}} : \psi &\edge b_i\circ \psi,
\qquad\qquad i=1,\dots, n-1, 
\\
T_{\al_i} : \psi &\edge t_i\circ\psi,
\qquad\qquad i=1,\dots, n,\\
Z: \psi &\edge z \circ\psi,
\end{aligned}
\end{equation}
compare with \reref{BTbt}.
\end{remark}
\subsection{}\lbb{subscf2}
Now we can prove that the complex $\Mt^\emptyset(\Si )$ is connected
and simply-connected. To prove that it is connected, it suffices to
check that the homeomorphisms $\psi, g\circ\psi$, where $g$ is one of
the generators of the group $\Ga_{0,n }$, can be connected by a path
in $\Mt^\emptyset(\Si )$.  This is obvious because $\Mt^\emptyset(\Si
)$ contains the edges~\eqref{btmoves}.

To prove that $\Mt^\emptyset(\Si )$ is simply-connected, note
first that it follows from \exref{xident3} that every path
can be deformed to a path that only uses $B_{\al, \be}$ for
neighboring boundary components $\al<\be$.
Therefore, any path can be contracted to a sequence of moves of the
form \eqref{btmoves}.  Thus, it remains to show that any closed loop
composed of the moves \eqref{btmoves} is contractible.  Since these
moves correspond to the generators of the group $\Ga_{0, n }$, used
in Proposition~\ref{pmcg0n}, this reduces to checking that the loops
corresponding to the relations (\ref{ga0n1}--\ref{ga0n8}) are
contractible.  This is straightforward.  The braid relation
\eqref{ga0n2} has already been established in \exref{xbrel}.  Using
\eqref{eqgenbr}, we can show that it suffices to check the relations
\eqref{ga0n7}, \eqref{ga0n8} for $n=2$, in which case they immediately
follow from the Dehn twist axiom. 
The other
relations follow from the commutativity of disjoint union
\eqref{comax} and the Cylinder axiom.  For example, both \eqref{ga0n3}
and \eqref{ga0n4} correspond to the identity $B_{\al_i, \al_{i+1}}
T_{\al_j} = T_{\al_j} B_{\al_i, \al_{i+1}}$.

This proves that $\Mt^\emptyset(\Si )$ is simply-connected, and
thus establishes assumption~2 of Proposition~\ref{pmpic}. 

\subsection{Finishing the proof}\lbb{subfinish}
So far, we have defined the map $\pi:\Mt(\Si)\to\C(\Si)$ and proved
that both the base and the fiber are connected and simply-connected,
thus establishing assumptions 1 and 2 of
Proposition~\ref{pmpic}. Assumption 4 is quite obvious, since the only
2-cells in $\C(\Si)$ are those obtained from the associativity axiom,
and they can be lifted to the 2-cells in $\Mt(\Si)$ also given by the
associativity axiom. Thus, the only thing that remains to be checked
is the assumption 3.

It is easy to see from the results of \ssref{subscf2} that any two
markings with the same cut system can be connected by a product of the
moves $Z, B_{\al,\be}$ ($\al,\be\in A(\Si_a)$ where $\Si_a$ 
is a connected component of $\Si\setminus C$), 
cf.\ \exref{xident3}.  

Thus, we only need to consider assumption 3 with $e_1$ being either
$Z$ or $B$. For $Z$, the statement immediately follows from the
symmetry of F axiom.

Hence, it suffices to check
that for $\Si=S_{0,n}\sqcup_c S_{0,k}$ and $\al,\be\in A(S_{0,n})$,
there exists
a path $e_2$ such that the following square is contractible in
$\Mt(\Si)$:
\begin{equation*}
\begin{CD}
M'_1 @>{F_c}>> M'_2 \\
@V{B_{\al,\be}}VV  @VV{e_2}V \\
M''_1 @>>{F_c}> M''_2
\end{CD}\quad.
\end{equation*}
This can be easily proved explicitly, using the axioms and
\exref{xident3}. Indeed, if both $\al$ and $\be$
are distinct from $c$, we can take $e_2=B_{\al, \be}$. If $\be=c$,
then we can take $e_2=B_{\al, I}$ where $I=A(S_{0,k})\setminus c$. 

Thus, we see that the map $\pi\colon\Mt(\Si)\to \C(\Si)$ satisfies all
assumptions of Proposition~\ref{pmpic} and thus, $\Mt(\Si)$ is
connected and simply-connected. This concludes the proof of \thref{tmtg0}.

\section{Proof of the main theorem for higher genus}\lbb{spfghi}
In this section we will prove the Main Theorem~\ref{tmtghi}
for higher genus: that for any extended surface $\Si$
the complex $\M(\Si)$, defined in Subsect.\ 
\ref{subg0mv}, \ref{subg0rel}, \ref{subghimv}, \ref{subghirel},
is connected and simply-connected.
The strategy of the proof is similar to the one used in the genus $0$ case.

First, we extend the complex $\M(\Si)$ by adding all 
disjoint unions and gluings of generalized
braidings and generalized S-moves as new edges, and adding their
definitions as new 2-cells. We denote this new complex by $\Mt(\Si)$;
again, $\M(\Si)$ is connected and simply-connected iff $\Mt(\Si)$ is
connected and simply-connected. Second, we define a complex $\Ct(\Si)$
with vertices the set $C(\Si)$ of all cut systems of $\Si$.  Then we
apply \prref{pmpic} to the canonical projection $\pi:\Mt(\Si)\to
\Ct(\Si)$. The most difficult part of the proof is checking that the
complex $\Ct(\Si)$ is simply-connected, which is based on the results
of \cite{HT} and \cite{H}.

\subsection{The complex $\C(\Si)$}\lbb{scsigma}
The definition of the complex $\Ct(\Si)$ is  parallel to the definition of  
$\Mt(\Si)$.
First, we define a complex $\C(\Si)$
with vertices the set $C(\Si)$ of all cut systems of $\Si$
(see \deref{dcut}).
The (directed) edges of $\C(\Si)$ are the following moves:

\begin{description}

\item[\barF-move] 
Let $\Si$ be an e-surface of genus zero, and let 
$C\in C(\Si)$ be a cut system on $\Si$, consisting of a single cut:
$C=\{c\}$.   Then we define  the
\barF-move $\bar F_c \colon C \edge \emptyset$, which removes $c$.

\item[\barS-move] Let $\Si$ be an e-surface of genus one with one
boundary component, and let $C$ be a cut system on $\Si$,
consisting of a single cut: $C=\{c\}$.  Let $c'$ be a simple closed
curve on $\Si$ such that $c'$ intersects $c$ transversally at exactly
one point (see \firef{Fsbarmv}). Then we add an edge $\bar S_{c,
c'}:\{c\}\edge \{c'\}$. 

\end{description}

\begin{figure}[h]
\begin{equation*}
\hspace{20pt}
\fig{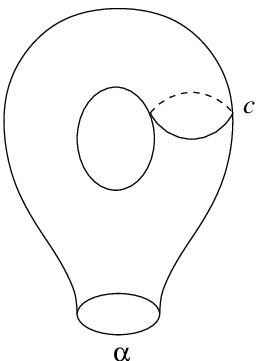}
\hspace{20pt}
\overset{\bar S}{\longedge}
\hspace{10pt}
\fig{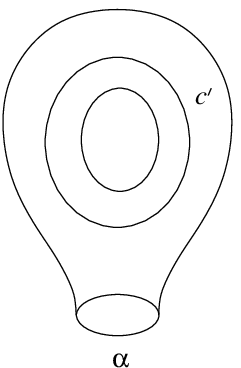}
\end{equation*}
\caption{\captionfont{\barS-move.}}\lbb{Fsbarmv}
\end{figure}

As before, we also add all the edges which can be obtained from the
\barF-, \barS-edges above by applying the operations of disjoint union
and gluing as in \ssref{sprop}. This implies that for every removable
cut $c\in C, C\in C(\Si)$ ($\Si$ not necessarily of genus zero), we
have an edge $\bar F_c:C\edge C\setminus\{c\}$. 

\begin{example}[{\bf Generalized \barS-move}]\lbb{xgensbarmv}
Let $\Si$ be a torus with $n$ holes, and let $c, c_1$ be cuts on $\Si$
as in \exref{xgensmv}.
 Then we  define the {\em generalized
$\bar S$-move\/}  as the composition of moves shown in \firef{Fgensmv} 
with $F,S$ replaced by $\bar F, \bar S$.
Again, it can be shown that the cut $c$ is uniquely
determined by $c_1, c_2$. Thus, we will denote this generalized
\barS-move by $\bar S_{c_1, c_2}$.  

It is easy to
see that if $C=\{c_1, \dots, c_k\}\in C(\Si)$, and $c_1'$ is a simple
closed curve on $\Si$ which intersects $c_1$ transversally at exactly
one point, and does not intersect any other cuts in $C$, then the
connected component $\Si_1$ of $\Si\setminus \{c_2, \dots, c_k\}$ which
contains $c_1, c'_1$ has genus one, and thus we have a generalized
\barS-move 
$\bar S_{c_1, c'_1}:\{c_1, c_2,\dots, c_k\}\edge \{c'_1, c_2,\dots,
c_k\}$ obtained by gluing the generalized  \barS-move on $\Si_1$ with the
identity on other components. 
\end{example}
\subsection{Relations in $\C(\Si)$}\lbb{srelincsi}
Let us impose the  following relations for the \barF- and \barS-moves:

\begin{description}
\item[Associativity of cut removal]
If $c_1,c_2\in C$ are two cuts on $\Si$
such that $\bar F_{c_1} \bar F_{c_2} (C)$ is defined, then
$\bar F_{c_2} \bar F_{c_1} (C)$ is defined and
\begin{equation}\lbb{barasax}
\bar F_{c_1} \bar F_{c_2} (C) = \bar F_{c_2} \bar F_{c_1} (C).
\end{equation}

\item[Inverse for \barS]
Let $\Si$ be a surface of genus one with one hole,
and $c$, $c'$ be as in the definition of \barS-move, cf.\
\firef{Fsbarmv}. Then
\begin{equation}\lbb{barsinv}
\bar S_{c', c} \bar S_{c, c'} (\{c\}) = \id(\{c\}).
\end{equation}

\item[Relation between \barS\ and \barF]
Let $\Si$ be a surface of genus one with two holes,
and $c_1$, $c_2$, $c_3$ be three cuts as shown in \firef{Fsfrel} below.
Then
\begin{equation}\lbb{sfrel}
\bar F_{c_1} \bar F_{c_2}^{-1} (\{c_1\}) 
= \bar S_{c_3, c_2} \bar S_{c_1, c_3} (\{c_1\}).
\end{equation}

\begin{figure}[h] 
\begin{equation*}
\fig{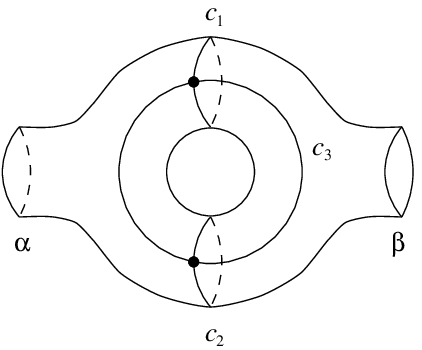}
\end{equation*}
\caption{\captionfont{Relation between \barS\ and \barF.}}\lbb{Fsfrel}
\end{figure}

\item[Triangle relation for \barS] 
Let $\Si$ be a torus with one hole,
and $c_1$, $c_2$, $c_3$ be three cuts as shown in
\firef{Fstriang} below.  Then
\begin{equation}\lbb{striang}
\bar S_{c_3, c_1} \bar S_{c_2, c_3} \bar S_{c_1, c_2} (\{c_1\}) 
= \id(\{c_1\}).
\end{equation}

\begin{figure}[h] 
\begin{equation*}
\fig{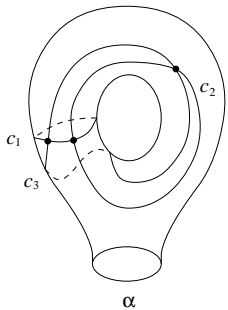}
\end{equation*}
\caption{\captionfont{Triangle relation  for \barS.}}\lbb{Fstriang}
\end{figure}

\item[Commutativity of disjoint union]
If $\Si=\Si_1\sqcup\Si_2$ and $E_i$ is an edge in $\C(\Si_i)$ ($i=1,2$),
then in $\C(\Si)$
\begin{equation}\lbb{comax2}
(E_1\sqcup\id)(\id\sqcup E_2)=(\id\sqcup E_2)(E_1\sqcup\id).
\end{equation}

\end{description}

Note that we also add the relations obtained from the above under
the action of the mapping class group;
for example, in Eq.\ \eqref{sfrel},
$c_1,c_2,c_3$ may be any three cuts such that 
$c_1,c_3$ and $c_2,c_3$ intersect at exactly one point and there are 
no other intersections.

Again, we add the propagation rules, i.e.\ we add all relations
obtained by taking disjoint unions and gluing,
cf.\ \ssref{subg0rel}.
Note that when $\Si$ is of genus $0$, the complex $\C(\Si)$
is the same as the one defined in \ssref{subcg0}.

Finally, as before, let us replace the complex $\C(\Si)$ by the
equivalent complex $\Ct(\Si)$, obtained by adding the generalized
\barS-moves as new edges (rather than considering them as
composition of moves), and adding the definition of these moves as new
relations.

\subsection{The projection}\lbb{pcproj}
We define the map of CW complexes $\pi\colon\Mt(\Si)\to\Ct(\Si)$, such
that on the vertices it is given by the canonical forgetting map
$\pi\colon M(\Si)\to C(\Si)$, and $\pi(Z)=\id$, $\pi(F) = \bar F$,
$\pi(B) = \id$, $\pi(S) = \bar S$.

Our goal is to prove that the projection map $\pi$ satisfies all the
assumptions of \prref{pmpic}. Obviously, as soon as we prove this, we
get a proof of the Main Theorem~\ref{tmtghi}. 

First of all, we need to check that the map $\pi$ is surjective on
edges, i.e.\ that every move in $\Ct(\Si)$ can be obtained by a projection of
a move in $\Mt(\Si)$. This is obvious for the \barF-move, and almost
obvious for the generalized \barS-move. 

\subsection{Checking assumption 2}\lbb{subass2}
Let us check that assumption 2 of \prref{pmpic} holds for the
projection map defined in \ssref{pcproj}. Clearly, for every cut
system $C$, $\pi^{-1}(C)\subset M(\Si)=\prod_a M^{\emptyset}(\Si_a)$
(cf.\ \ssref{subscfiber}), and a path $l$ which lies in
$\pi^{-1}(C)$ must be composed of Z-, B-moves only 
(in particular, it
cannot include an S-move). Thus, the same proof as in the genus zero case
(see Subsect.~\ref{subscfiber}--\ref{subscf2}) applies here.

\subsection{Checking assumption 3}\lbb{subass3}
Let us check that assumption 3 of \prref{pmpic} holds for the
projection map defined in \ssref{pcproj}, i.e.\ that for every
edge $e:C_1\edge C_2$ and  two  its liftings $e', e''$ to $\Mt$, they
can be included in a commutative square. If the edge $e$ is of the
\barF -type, then the same proof as in the genus $0$ case (see
\ssref{subfinish})  applies. 

Thus, we have to consider the case when the
edge $e$ is of \barS -type. This reduces to asking what different
liftings a given generalized \barS-move has.  This is answered by
the following lemma.
\begin{lemma} \lbb{lslift}
Let $\Si$ be a torus with $n$ holes $\al_1,\dots,\al_n$,
and let $c$ be a cut on it as shown
in \firef{FstmarkforS}. Let $M=(\{c\},m)$ be a marking on $\Si$ such
that $S(M)=(\{c'\}, m')$, with the cut $c'$ shown in
\firef{FstmarkforS} {\rm(}recall that the generalized $S$-move is uniquely
defined by $M${\rm)}. Then any such $M$ can be connected by a sequence of
moves $B_{\al_i, \al_{i+1}}$ $(i=1,\dots, n-1)$ and
$T_{\al_i}$ $(i=1,\dots, n)$ and their inverses with one of the
two standard markings shown in \firef{FstmarkforS}.
\end{lemma}
\begin{figure}[h]
\begin{equation*}
\fig{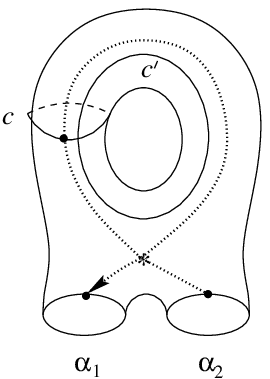}
\hspace{60pt}
\fig{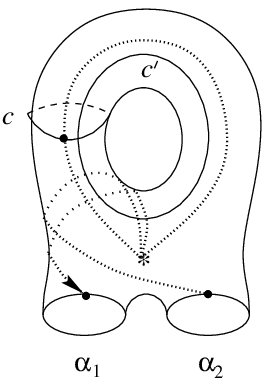}
\end{equation*}
\caption{\captionfont{Two standard markings of a torus with $n$ holes 
(for $n=2$).}}
\lbb{FstmarkforS}
\end{figure}

\begin{proof} Denote the two markings in \firef{FstmarkforS} by $M',
  M''$. Now, let $M$ be the marking satisfying the conditions of the
  theorem. Then there exists a homeomorphism $\ph':\Si\to \Si$ such
  that $\ph'(M)=M'$ and $\ph' (c)=c$. Moreover, it is easy to see
  that we must also have $\ph'(c')=c'$.

Presenting the torus as a rectangle with identified opposite sides, we
see that such a $\ph'$ is the same as a homeomorphism of a rectangle
with $n$ holes onto itself which maps vertical sides to vertical and
horizontal to horizontal. If $\ph'$ preserves each of the sides, then
without loss of generality we may assume that it acts as identity on
the boundary of the rectangle. But it is well-known that the group of 
such homeomorphisms is generated by the elements
$b_i$ $(i=1,\dots, n-1)$ and  $t_j$ $(j=1,\dots, n)$, 
cf. \cite[Theorem 1.10]{Bir}. Thus, in this case $M$ and $M'$ can be
connected by a sequence of $B, T$ moves as in the theorem. 

If the homeomorphism $\ph'$ interchanges the opposite sides of
rectangle (i.e., \ interchanges the sides of the cuts $c,c'$), then
we need to repeat the same argument for $M''$; it is easy to check
that in this case the homeomorphism  $\ph''$ will preserve each of the
sides of the rectangle. 
\end{proof}

Arguing as in \ssref{subfinish}, we
see that it suffices to check that we can find a path $e_2$ which does not
change the cut system and such that $Se_1=e_2S$, with $e_1$ being
$B_{\al_i, \al_{i+1}}$ $(i=1,\dots, n-1)$ or
$T_{\al_i}$ $(i=1,\dots, n)$. This is
obvious, since we can take $e_2=e_1$, and the equality would follow
from the commutativity of disjoint union. And, finally, it remains to
show that we can find $e_1, e_2$ such that the square 
\begin{equation*}
\begin{CD}
M' @>{S}>>  M_2'  \\
@V{e_1}VV  @VV{e_2}V \\
M'' @>{S}>>  M_2'' 
\end{CD}
\end{equation*}
is commutative,
where $M'_1, M''_1$ are the standard markings in
\firef{FstmarkforS}. But this can be easily achieved by letting
$e_1=Z^{-1}B_{\{\al_1, \dots, \al_n\}, c}$, $e_2=Z^{-1}B_{\{\al_1, \dots,
\al_n\}, c'}$.  Indeed, it follows from \eqref{g1n1r1} that $e_1=S^2,
e_2=S^2$ (cf.\ \firef{Fg1n1rel}), and thus the square above is
obviously commutative. 

\subsection{Checking assumption 4}\lbb{subass4}
Let us check that assumption 4 of \prref{pmpic} holds for the
projection map defined in \ssref{pcproj}, i.e.\ that for every
every 2-cell $X$ in $\Ct(\Si)$, its boundary can be lifted to a contractible
loop in $\Mt(\Si)$. In other words, we need to check that every relation in
$\Ct(\Si)$ can be obtained by projecting some relation in $\Mt(\Si)$. 
Clearly, it suffices to check this for  the basic relations
(\ref{barasax}--\ref{striang}). 

For the associativity axiom \eqref{barasax}, this is obvious: it can
be obtained by projecting the associativity axiom \eqref{asax1}. The
inverse axiom \eqref{barsinv} can be easily obtained from 
the relation $S^2=Z^{-1}B$ (see \eqref{g1n1r1} and \firef{Fg1n1rel}). 
The relation \eqref{sfrel} between $\bar S$ and $\bar F$ 
is nothing else but the projection of the defining relation
\eqref{g1n2}, see Appendix~\ref{apg1n2}.
Similarly, the triangle
relation \eqref{striang} is exactly  the projection of the
relation $(ST)^3=S^2$ (see \eqref{g1n1r2} and Appendix~\ref{apg1n1r2}). 

\begin{theorem}[$g\ge0$]\lbb{tcsighi}
The complex $\Ct(\Si)$ is connected and simply-connected.
\end{theorem}
This theorem is proved in Subsections~\ref{pcmin}, \ref{spftcsighi} below.

Without loss of generality we may assume that $\Si$
is connected.
Recall that a cut system is called minimal if it contains no removable
cuts; this is exactly what is called a ``cut system'' in \cite{HT, H}.
Let $\Ct^\Min(\Si)$ be the complex with vertices: all minimal cut
systems, edges: all generalized \barS-moves, and the relations
induced by the relations in $\Ct(\Si)$ (i.e., a path in $\Ct^\Min(\Si)$
is contractible if it is contractible as a path in $\Ct(\Si)$).

\begin{proposition}\lbb{pcmin}
 The subcomplex $\Ct^\Min(\Si)$ is connected and
simply-connected.
\end{proposition}
\begin{proof} 
The proof is based on the results of \cite[Section~2]{H}, where 
a certain $2$-dimensional CW complex $Y_2$ is introduced,
which has the same vertices and
edges as $\Ct^\Min(\Si)$, but different 2-cells. 
Since $Y_2$ is connected and
simply-connected \cite[Theorem~2.2]{H}, 
it suffices to show that all the relations of 
Harer  follow from the relations in $\Ct(\Si)$. 

The first relation \cite[Eq.\ (R${}_1$)]{H} has the form 
\begin{equation}\lbb{harer1}
\bar S_{c_3, c_1} \bar S_{c_2, c_3} \bar S_{c_1, c_2} (\{c_1\})
= \id(\{c_1\}),
\end{equation}
where $c_1, c_2, c_3$ are some cycles on a surface $\Si$.  There are
many different choices, displayed in \cite[Figure~4]{H}.  However, it
is easy to see that, by making one additional cut, they all reduce to
the configuration shown in \firef{F1harera}.

\begin{figure}[h]
\begin{equation*}
\fig{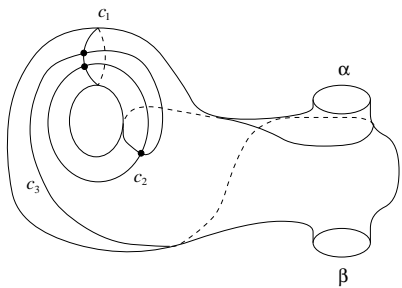}
\end{equation*}
\caption{\captionfont{Harer's first relation.}}\lbb{F1harera}
\end{figure}

To prove \eqref{harer1} for the cuts shown in \firef{F1harera}, 
redraw \firef{F1harera} as shown in
\firef{F1harerb}, and add one more cut $c_4$.

\begin{figure}[h]
\begin{equation*}
\fig{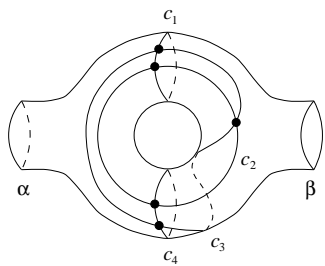}
\end{equation*}
\caption{\captionfont{Proof of the Harer's first relation.}}\lbb{F1harerb}
\end{figure}

Then consider the following diagram:
\begin{equation*}
\xymatrix{
& \{1\} \ar[d]\ar[ddr]^{\bar S}&\\
& \{4\} \ar[dl]^{\bar S} \ar[dr]_{\bar S}&\\
\{3\}\ar[uur]^{\bar S} & &\ar[ll]^{\bar S}\{2\}
}
\end{equation*}
where $\{1\}$ stands for the cut system consisting of one cut $c_1$,
etc., and the vertical arrow $\{1\}\to \{4\}$ is given by
$\bar F_{c_1}\bar F^{-1}_{c_4}$. 
 The outer triangle is exactly the left-hand side of the
relation \eqref{harer1}. 
On the other hand, the two top small triangles are
contractible by \eqref{sfrel}, and the bottom triangle is contractible
because it can be obtained from the Triangle relation \eqref{striang} by
gluing a three-punctured sphere to the hole. This completes the proof
of the first Harer's relation \eqref{harer1}.

The second relation \cite[Eq.\ (R${}_2$)]{H}
states that if $c_1,\dots,c_4$ are $4$ cycles such that $c_1$
intersects $c_2$ at one point, $c_3$ intersects $c_4$ at one point,
and there are no other intersections (this is illustrated by the
diagram in \firef{F2harer}), then
\begin{equation*}
 \bar S_{12} \bar S_{34} = \bar S_{34} \bar S_{12},
\end{equation*}
where $\bar S_{ij} = \bar S_{c_i, c_j}$. 
This follows from the commutativity of disjoint union \eqref{comax2}.

\begin{figure}[h]
\begin{equation*}
\fig{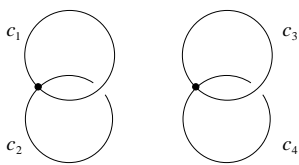}
\end{equation*}
\caption{\captionfont{Harer's second relation.}}\lbb{F2harer}
\end{figure}

The third relation \cite[Eq.\ (R${}_3$)]{H} is
\begin{equation*}
\bar S_{42} \bar S_{54} \bar S_{61} 
\bar S_{26} \bar S_{35} \bar S_{13} (\{c_1,c_2\})= \id(\{c_1,c_2\}),
\end{equation*}
where $c_1,\dots,c_6$ are cuts on an e-surface $\Si$ of genus two
with one hole, displayed in \firef{F3harer} (cf.\ \cite[Figure~3]{H}).

\begin{figure}[h]
\begin{equation*}
\fig{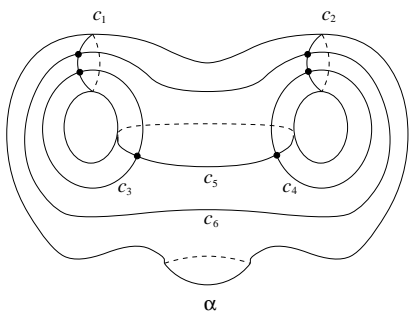}
\end{equation*}
\caption{\captionfont{Harer's third relation.}}\lbb{F3harer}
\end{figure}

This relation follows from \eqref{sfrel}:
\begin{equation*}
\spreaddiagramcolumns{-1pc}\spreaddiagramrows{0pc}
{\diagram
& \{1,2\} \rrto^-{\bar S} \ar@<1ex>[dr]^-{\bar F^{-1}} & &
\{2,3\} \drto^-{\bar S} &\\
\{1,4\} \urto^-{\bar S} & & 
\quad \{1,2,5\} \ar@<1ex>[ul]^-{\bar F} \ar@<1ex>[dl]^-{\bar F} 
\ar@<1ex>[rr]^-{\bar F} & &
\{2,5\} \dlto^-{\bar S} \ar@<1ex>[ll]^-{\bar F^{-1}}  \\ 
& \{1,5\} \ulto^-{\bar S} \ar@<1ex>[ur]^-{\bar F^{-1}} & &
\{5,6\} \llto^-{\bar S} &\\
\enddiagram}
\quad .
\end{equation*}

This completes the proof that $\Ct^\Min(\Si)$ is simply-connected.
\end{proof}
\subsection{Proof of \thref{tcsighi}}\lbb{spftcsighi}
Since any cut system can be joined to a minimal one by 
erasing cuts, and $\Ct^\Min(\Si)$ is connected, 
it follows that $\Ct(\Si)$ is connected.

To prove that $\Ct(\Si)$ is simply-connected, we first note that every path
$\bar F_{c_1} \bar F_{c_2}^{-1}$ is homotopic to either
$\id$, or $\bar F_{c_2}^{-1} \bar F_{c_1}$, or 
$\bar S_{c_3,c_2} \bar S_{c_1,c_3}$ for certain $c_3$. 
Indeed, if both $\bar F_{c_1} \bar F_{c_2}^{-1} (C)$ and
$\bar F_{c_2}^{-1} \bar F_{c_1} (C)$ are defined, then they are
equal. Suppose that the first one is defined but the second one is not.
Then $c_1$ and $c_2$ do not intersect and $\Si\setminus C$ 
becomes of positive genus
if we remove them from the cut system $C$.
Hence, there is a cut $c_3$ which intersects both of them,
and we can apply \eqref{sfrel}.

Similarly, note that
every path $\bar S_{c_2,c_3} \bar F_{c_1}^{-1}$
is homotopic to
$\bar F_{c_1}^{-1} \bar S_{c_2,c_3}$.
Indeed, 
$\bar S_{c_2,c_3} \bar F_{c_1}^{-1} (C)$ 
being defined implies that $c_1$ intersects
neither $c_2$ nor $c_3$. Then by the commutativity of disjoint union,
$\bar F_{c_1}^{-1} \bar S_{c_2,c_3}(C)$ is also defined and they are
equal.

Now take any closed loop $l$ in $\Ct(\Si)$. Without loss of generality,
we may assume that it has a minimal cut system as the basepoint.
Using the above two remarks, we can deform $l$ into a loop 
composed only of \barS-moves.
Indeed, we can move any \barF\inv-move to the left until it
meets an \barF-move and either cancels out or creates a pair of \barS-moves.
Repeating this procedure, 
we will get a loop composed only of \barF- and \barS-moves.
But since the number of cuts in the initial cut system should be the same,
it is actually composed only of \barS-moves.

If we start with a minimal cut system and apply to it a sequence of
\barS-moves, we again get a minimal cut system. Since the subcomplex
$\Ct^\Min(\Si)$ of $\Ct(\Si)$ is simply-connected by \prref{pcmin}, it
follows that $l$ is contractible.

This completes the proof that $\Ct(\Si)$ is
simply-connected. Therefore, we have checked all the assumptions of
\prref{pmpic}, and thus, we have proved that the complex $\M(\Si)$ is
connected and simply-connected. 

\appendix

\section{The relation \eqref{g1n1r2}: $TSTST=S$}\lbb{apg1n1r2}

\begin{equation*} 
\hspace{-20pt}
\fig{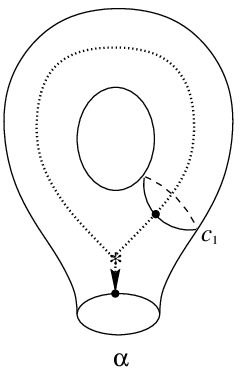}
\hspace{10pt}
\overset{T}{\longedge}
\hspace{10pt}
\fig{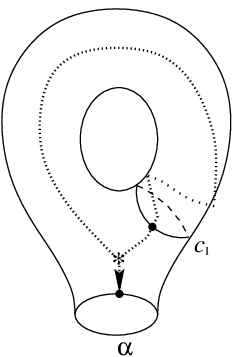}
\hspace{10pt}
\overset{S}{\longedge}
\hspace{10pt}
\fig{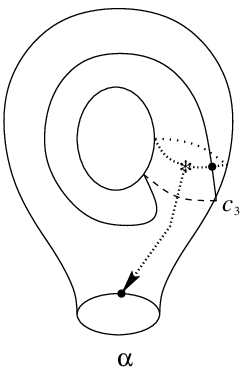}
\hspace{10pt}
\overset{T}{\longedge}
\end{equation*}

\begin{equation*} 
\overset{T}{\longedge}
\hspace{5pt}
\fig{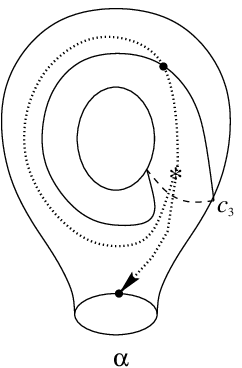}
\hspace{10pt}
\overset{S}{\longedge}
\hspace{10pt}
\fig{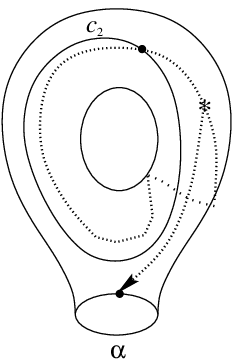}
\hspace{10pt}
\overset{T}{\longedge}
\hspace{10pt}
\fig{smv2.eps}
\end{equation*}

\section{The relation \eqref{g1n2}:
$B_{\al,\be} F_{c_1} F_{c_2}^{-1} = S^{-1} {\wti T}^{-1} T S$}\lbb{apg1n2}

The left hand side of \eqref{g1n2} is:

\begin{align*} 
\fig{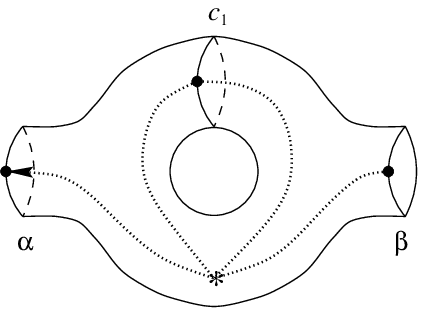}
\hspace{10pt}
&\overset{F_{c_2}^{-1}}{\longedge}
\hspace{10pt}
\fig{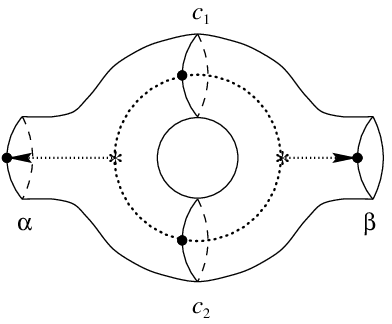}
\hspace{10pt}
\overset{F_{c_1}}{\longedge}
\\
\overset{F_{c_1}}{\longedge}
\hspace{10pt}
\fig{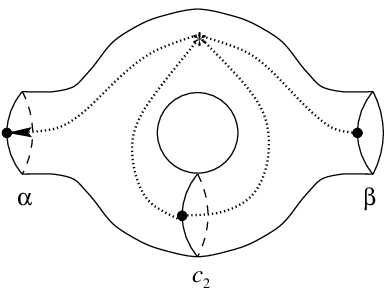}
\hspace{10pt}
&\overset{B_{\al,\be}}{\longedge}
\hspace{10pt}
\fig{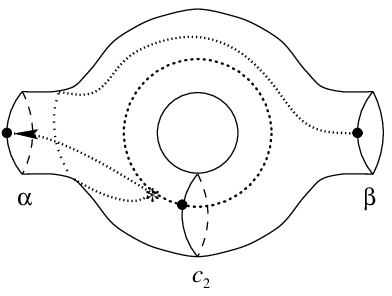}
\end{align*}

The right hand side of \eqref{g1n2} is:

\begin{gather*} 
\fig{b2.eps}
\hspace{10pt}
\overset{S}{\longedge}
\hspace{10pt}
\fig{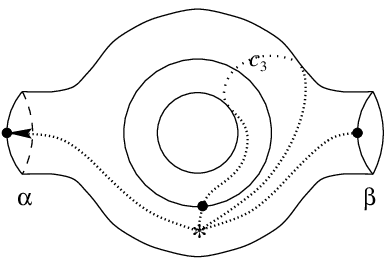}
\hspace{10pt}
\overset{T_{c_3+\be}^{-1} T_{c_3}}{\longedge}
\\
\overset{T_{c_3+\be}^{-1} T_{c_3}}{\longedge}
\hspace{10pt}
\fig{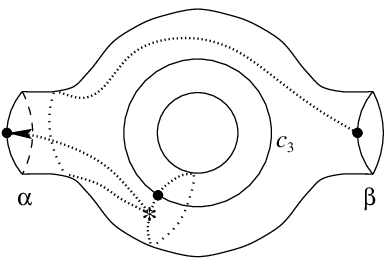}
\hspace{10pt}
\overset{S^{-1}}{\longedge}
\hspace{10pt}
\fig{b9.eps}
\end{gather*}

Below is a version of the same relation which makes sense in the
complex $\M^\Max$, i.e.\ which only uses spheres with $\le 3$ 
holes---for brevity, we just wrote the corresponding marking graphs.

\begin{equation*} 
\fig{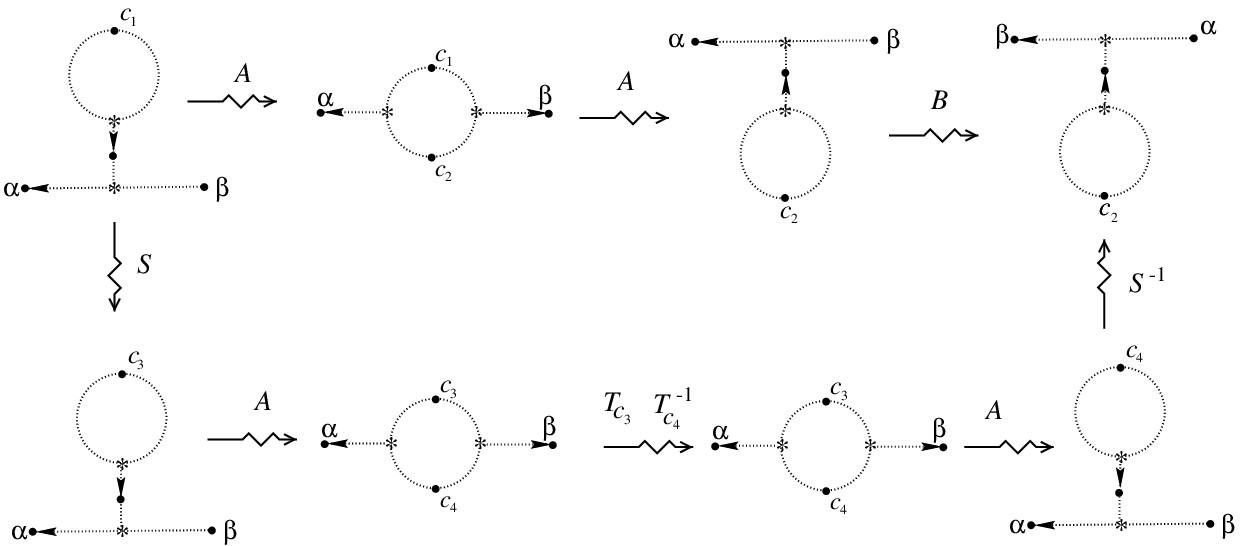}
\end{equation*}

\section{Triangle, Pentagon, and Hexagon relations}\lbb{appenthex}

In this appendix, we formulate the Triangle, Pentagon, and Hexagon
axioms, which were used in Section~\ref{smmaxsi}. For brevity, we only
give pictures of the corresponding marking graphs; as was mentioned
before, this is sufficient to uniquely reconstruct the moves.  All
unmarked edges in these diagrams are compositions of the form $
(Z^*\sqcup Z^*) A (Z^*\sqcup Z^*)$; the powers of $Z$ are uniquely
determined by the distinguished edges in the diagram and by the
requirement that this composition is well-defined.

The Triangle axiom requires that the 
diagram in \firef{Ftriangle} below be commutative.

The Pentagon relation is shown in \firef{Fpent}.

Finally, there are two Hexagon axioms. One of them claims the
commutativity of the diagram in \firef{Fhexagon2}; 
the other is obtained by replacing
all occurrences of $B$ by $B^{-1}$, so that $B_{\al\be}$ is replaced by
$B_{\be\al}^{-1}$, etc. 

\begin{figure}[h]
\begin{equation*}
\fig{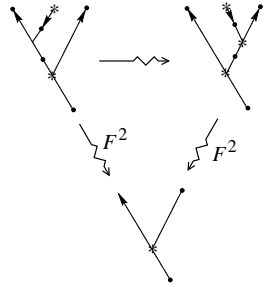}
\end{equation*}
\caption{\captionfont{Triangle relation.}}\lbb{Ftriangle}
\end{figure}

\begin{figure}[h]
\begin{equation*}
\fig{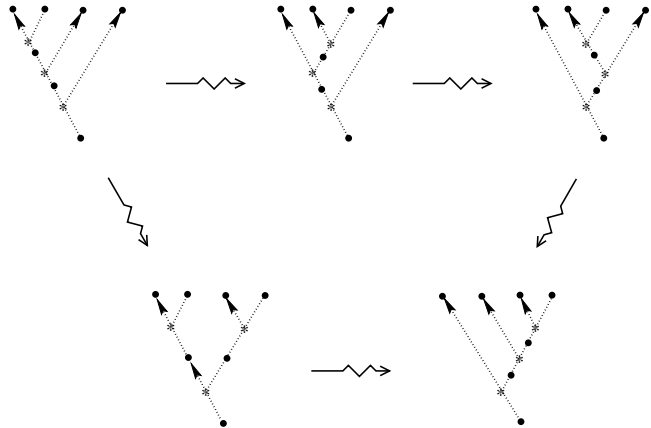}
\end{equation*}
\caption{\captionfont{Pentagon relation.}}\lbb{Fpent}
\end{figure}

\begin{figure}[h]
\fig{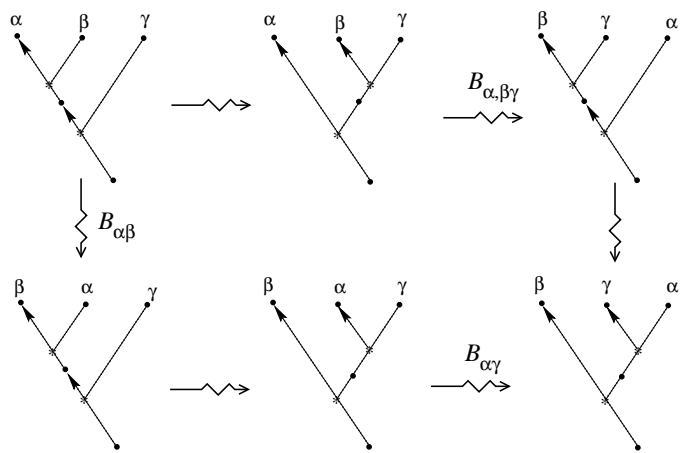}
\caption{\captionfont{Hexagon relation.}}\lbb{Fhexagon2}
\end{figure}

\bibliographystyle{amsalpha}

\end{document}